\documentclass[12pt,reqno,tikz]{amsart}
\usepackage{amsopn,amsfonts,amssymb,amsthm,mathtools,amsmath}
\usepackage{tikz}
\usepackage{tikz-cd}
\usetikzlibrary{trees,arrows}
\usepackage[edges]{forest}
\usepackage{latexsym}

\usepackage[T1]{fontenc}
\usepackage{lmodern}
\usepackage{textcomp}
\usepackage[hidelinks]{hyperref}
\usepackage[utf8]{inputenc}
\usepackage{enumerate}
\usepackage[shortlabels]{enumitem}
\setlist[enumerate, 1]{\sc(1)}
\usepackage{float}
\usepackage{hhline}
\usepackage{multirow}
\usepackage{anysize}

\usepackage{enumerate}

\def\z{\mathfrak{z}}
\def\u{\mathfrak{u}}
\def\b{\mathfrak{b}}

\def\g{\mathfrak{g}}

\def\n{\mathfrak{n}}

\def\a{\mathfrak{a}}

\def\so{\mathfrak{so}}

\def\u{\mathfrak{u}}
\def\gl{\mathfrak{gl}}

\def\M{\mathcal{M}}

\def\B{\mathcal{B}}

\def\l{\ell}

\def\C{\mathbb{C}}
\def\R{\mathbb{R}}
\def\Q{\mathbb{Q}}
\def\Z{\mathbb{Z}}
\def\N{\mathbb{N}}

\newcommand{\ort}{{\sf O}}
\newcommand{\GL}{{\sf GL}}
\newcommand{\SL}{{\sf SL}}
\newcommand{\SU}{{\sf SU}}

\def\ad{\operatorname{ad}}

\def\tr{\operatorname{tr}}

\def\I{\operatorname{I}}
\def\Id{\operatorname{Id}}

\newcommand{\Aut}{\operatorname{Aut}}
\newcommand{\Coker}{\operatorname{Coker}}
\newcommand{\Ker}{\operatorname{Ker}}

\DeclareMathOperator{\Iso}{Iso}
\DeclareMathOperator{\im}{Im}
\DeclareMathOperator{\ord}{ord}
\DeclareMathOperator{\hol}{Hol}

\def\alt{\raise1pt\hbox{$\bigwedge$}}

\def\la{\langle}
\def\ra{\rangle}
\def\multiset#1#2{\ensuremath{\left(\kern-.3em\left(\genfrac{}{}{0pt}{}{#1}{#2}\right)\kern-.3em\right)}}

\DeclareSymbolFont{extraup}{U}{zavm}{m}{n}

\theoremstyle{plain}
\newtheorem{teo}{\bf Theorem}[section]
\newtheorem{coro}[teo]{\bf Corollary}
\newtheorem{prop}[teo]{\bf Proposition}
\newtheorem{lema}[teo]{\bf Lemma}

\theoremstyle{definition}
\newtheorem{defi}[teo]{\bf Definition}
\newtheorem{ejem}[teo]{\bf Example}
\theoremstyle{remark}
\newtheorem{obsx}[teo]{Remark}
\newenvironment{obs}
{\pushQED{\qed}\obsx}

\newcommand{\ri}{{\rm (i)}}
\newcommand{\rii}{{\rm (ii)}}
\newcommand{\riii}{{\rm (iii)}}

\newcommand{\matriz}[1]{\ensuremath{\begin{pmatrix}#1\end{pmatrix}}}

\newcommand{\comillas}[1]{``#1''}

\title{Classification of 6-dimensional splittable flat solvmanifolds}
\author{Alejandro Tolcachier}
\email{atolcachier@famaf.unc.edu.ar}

\address{FaMAF, Universidad Nacional de C\'ordoba, Ciudad Universitaria, X5000HUA C\'ordoba, Argentina}
\begin{document}
	
	\renewcommand{\bibname}{References}
	
	\begin{abstract}
		A flat solvmanifold is a compact quotient $\Gamma\backslash G$ where $G$ is a simply-connected solvable Lie group endowed with a flat left invariant metric and $\Gamma$ is a lattice of $G$. Any such Lie group can be written as $G=\R^k\ltimes_{\phi} \R^m$ with $\R^m$ the nilradical. In this article we focus on 6-dimensional splittable flat solvmanifolds, which are obtained quotienting $G$ by a lattice $\Gamma$ that can be decomposed as $\Gamma=\Gamma_1\ltimes_{\phi}\Gamma_2$, where $\Gamma_1$ and $\Gamma_2$ are lattices of $\R^k$ and $\R^m$, respectively. We analyze the relation between these lattices and the conjugacy classes of finite abelian subgroups of $\GL(n,\Z)$, which is known up to $n\leq 6$. From this we obtain the classification of 6-dimensional splittable flat solvmanifolds.

	\end{abstract}
	
	\maketitle

	
	\section{Introduction}
	
	A solvmanifold is a compact quotient $\Gamma\backslash G$ of a simply-connected solvable Lie group $G$ by a discrete subgroup $\Gamma$ (such a cocompact discrete subgroup $\Gamma$ is called a lattice). When this solvmanifold carries a flat Riemannian metric (i.e., its Levi-Civita connection has curvature zero) induced by a flat left invariant Riemannian metric on $G$ then $\Gamma\backslash G$ is called a flat solvmanifold. In this case $\Gamma\backslash G$ is a compact flat manifold and since $\Gamma=\pi_1(\Gamma\backslash G)$, $\Gamma$ is a Bieberbach group.
	
	Questions concerning the classification of flat manifolds and of solvmanifolds have been studied since the early 20th century (see the book of Charlap \cite{Ch} on flat manifolds, and for instance Chevalley \cite{Che}, Malcev \cite{Malcev}, Mostow \cite{Mostow1}, Auslander \cite{Aus601}, on solvmanifolds).
	
	More recently, compact solvmanifolds have been used as a source of examples and counterexamples in diverse areas of differential geometry. Solvmanifolds generalize the well known family of
	nilmanifolds which are defined similarly when $G$ is nilpotent. Many important global properties, such as cohomological properties, of nilmanifolds cannot be generalized to solvmanifolds. For instance, the well known Nomizu's theorem \cite{Nomizu} which allows to compute the de Rham cohomology of a nilmanifold in terms of the cohomology of the Lie algebra does not necessarily hold for a solvmanifold. 
	Moreover, it is difficult to determine whether a given solvable Lie group admits lattices or not. Nevertheless, there is such a criterion for nilpotent Lie groups and it was given by Malcev in \cite{Malcev}.
	
	Solvmanifolds and nilmanifolds have also had several applications in complex geometry. For instance, the Kodaira-Thurston nilmanifold was the first example of a compact non-Kähler symplectic manifold \cite{KodTh}. The completely-solvable Nakamura manifold is an example of a cohomologically Kähler non-Kähler solvmanifold \cite{FMS}. The well known Oeljeklaus-Toma manifolds (which are compact complex non-Kähler manifolds constructed by using number theory in \cite{OT}) were described by Kasuya as solvmanifolds and using this description he proved that they do not admit any Vaisman metric \cite{Kasuya}.
	
	On the other hand, compact flat manifolds are well understood due to the classical Bieberbach's theorems and they have been used to study different phenomena in geometry. For instance, questions about isospectrality (see \cite{MiaRo} and the references therein), Kähler flat metrics with holonomy in $\SU(n)$ (\cite{Dekimpe}), among others. 
	
	The class of flat solvmanifolds lies in the intersection between the two well studied theories of solvmanifolds and compact flat manifolds, and thus provide a nice interplay between them. Also, this class is rich enough to produce a diverse collection of examples.  
	
	In \cite{AusAus}, L. Auslander and M. Auslander characterized the fundamental groups of compact flat solvmanifolds. In \cite{Morgan}, A. Morgan gave a classification of such manifolds (in the more general case when $\Gamma$ is closed but not necessarily discrete) in dimensions $n\leq 5$, except the
	five-dimensional ones with first Betti number equal to one.
	
	In a previous paper (\cite{Tol}), we studied flat solvmanifolds from the point of view of the well known Milnor's characterization of solvable Lie groups which admit a flat left-invariant metric. We combined this with Bieberbach's classical theory to prove some properties of the holonomy group of a flat solvmanifold. Moreover, we also gave a classification of flat solvmanifolds in dimensions 3, 4 and 5, using the fact that in these dimensions all these solvmanifolds are quotients of almost abelian Lie groups,  i.e., simply-connected solvable Lie groups whose Lie algebra has a codimension-one abelian ideal. 
	
	The main goal of this article is to continue the study initiated in \cite{Tol} and, in particular, to provide the classification of splittable flat solvmanifolds in dimension 6. In \S3 we outline some facts about splittable solvmanifolds. A solvmanifold $\Gamma\backslash G$ will be called \textit{splittable} if $G$ is a semidirect product $G=\R^k\ltimes_{\phi} N$, where $N$ is the nilradical of $G$, and the lattice $\Gamma$ is splittable, i.e., $\Gamma=\Gamma_1\ltimes_{\phi} \Gamma_2$ where $\Gamma_1$ and $\Gamma_2$ are lattices in $\R^k$ and $N$ respectively. A flat Lie group $G$ is splittable since it can be written, according to \cite{Mi}, as $G=\R^k\ltimes_{\phi} \R^m$ where $\R^m$ is the nilradical of $G$. 
	
	After recalling a criterion given by \cite{Y} to determine all the splittable lattices in a splittable Lie group, we embed the holonomy group of a $(k+m)$-dimensional splittable flat solvmanifold as a subgroup of $\GL(m,\Z)$. From this point of view we can give a relation between the splittable lattices of a flat Lie group and finite abelian subgroups of $\GL(k,\Z)$, which allows us to give a way to classify splittable flat solvmanifolds. We point out that the classification of finite subgroups of $\GL(n,\Z)$ is known up to $n\leq 6$.
	
	In $\S4$ we study the 6-dimensional case. In dimension 6 any flat Lie group is of the form $G=\R\ltimes \R^5$ or $G=\R^2\ltimes \R^4$. In \S4.1 we obtain a complete classification of 6-dimensional almost abelian solvmanifolds. In this case, we study the conjugacy classes of matrices of finite order in $\GL(5,\Z)$ and we obtain the conjugacy classes of subgroups. We do this by applying the methods developed in \cite{Yang}, without relying on any known classification, to illustrate a possible way to tackle the problem in higher dimensions. Finally, in \S4.2 we obtain a complete classification of 6-dimensional splittable non almost abelian solvmanifolds, by using the classification of conjugacy classes of finite subgroups of $\GL(4,\Z)$.
	
	\medskip
	
	\textit{Acknowledgements.} I am very grateful to my Ph.D. advisor Adrián Andrada for his continuous guidance during the writing of this paper. I also thank Jonas Deré and Derek Holt for the useful conversations and Examples 3.3 and 3.8, respectively.
	
	\section{Preliminaries on compact flat manifolds and solvmanifolds}
	
	A discrete and torsionfree subgroup $\Gamma$ of $\Iso(\R^n)\cong \ort_n\ltimes \R^n$ with the property that the orbit space $\R^n/\Gamma$ is compact is called a \textit{Bieberbach} group.
	
	An equivalent condition for $\Gamma$ to be torsionfree is that the action of $\Gamma$ is free. Moreover, if a subgroup of $\Iso(\R^n)$ is discrete and acts freely on $\R^n$ then the action is also properly discontinuous (for all $x\in\R^n$ there exists an open neighborhood $U_x$ such that $g U_x\cap U_x=\emptyset$ for all $g\neq \Id$) and the converse also holds. In this case, $\R^n/\Gamma$ admits the structure of a differentiable manifold and $\pi_1(\R^n/\Gamma)\cong\Gamma$.
	
	An important subgroup of a Bieberbach group $\Gamma$ is the subgroup $\Lambda$ of pure translations of $\Gamma$, i.e., the elements $(A,v)\in \Gamma$ such that $A=\I$. Note that $\Lambda=\Gamma\cap \R^n$.
	
	Bieberbach groups are well described by three classical theorems known as \comillas{Bieberbach's theorems}. For more details about Bieberbach groups, see \cite{Ch}.
	
	\begin{teo}[Bieberbach's First Theorem]\label{maximal}
		Let $\Gamma$ be a Bieberbach subgroup of $\Iso(\R^n)$ and $\Lambda=\Gamma\cap \R^n$. Then $\Lambda$ is a normal free abelian subgroup of rank $n$ and $\Gamma/\Lambda$ is a finite group. Furthermore, $\Lambda$ is the unique maximal normal abelian subgroup.
	\end{teo}
In other words, a Bieberbach group $\Gamma$ satisfies an exact sequence $1\to \Lambda\xrightarrow[]{\iota} \Gamma\xrightarrow[]{\pi} H\to 1$, where $\Lambda$ is the traslation subgroup of $\Gamma$ and $H=\Gamma/\Lambda$ is a finite group, called the \textit{holonomy group} (or \textit{point group}) of $\Gamma$, which can be identified to $p_1(\Gamma)$, where $p_1:\Gamma\to \ort(n)$ is the projection into the first factor. Since $\Lambda$ is abelian, the action $H$ on $\Gamma$ is given by $h\cdot \gamma=\tilde{h} \gamma \tilde{h}^{-1}$, where $\tilde{h}$ is any element of $\Gamma$ satisfying $\pi(\tilde{h})=h$.

Following \cite{Hi1986}, the set of data $(H,\Lambda):=(H,\Lambda,\cdot)$ is called a \textit{crystal class}. More precisely, $(H,\Lambda)$ is the set of all $n$-dimensional Bieberbach groups $\Gamma$ that appear as an extension of $H$ by $\Lambda$. Two crystal classes $(H,
\Lambda)$ and $(H',\Lambda')$ are \textit{arithmetically equivalent} if there are isomorphisms $\varphi:H\to H'$ and $\alpha:\Lambda\to\Lambda'$ so that 
\begin{equation}\label{alpha} \alpha(h\cdot\gamma)=\varphi(h)\cdot \alpha(\gamma)\, \text{for all}\, \gamma\in \Lambda.
	\end{equation} Choosing a free integral basis for $\Lambda$, a crystal class can be viewed as a one-to-one homomorphism $H\to \Aut(\Lambda)\cong \GL(n,\Z)$. So after choosing such bases for $\Lambda$ and $\Lambda'$, $H$ and $H'$ can be considered as  subgroups of $\GL(n,\Z)$ and the condition \eqref{alpha} transforms into $\alpha H \alpha^{-1}=H'$, where $\alpha\in \GL(n,\Z)$. Therefore $H$ and $H'$ are conjugated in $\GL(n,\Z)$. The resulting equivalence classes are the \textit{arithmetic crystal classes} or \textit{$\Z$-classes}.

\begin{obs}\label{biebiso}
	Isomorphic Bieberbach groups determine the same arithmetic crystal class. Indeed, if $F:\Gamma_1\to \Gamma_2$ is an isomorphism, then $F(\Lambda_1)=\Lambda_2$ where $\Lambda_i$ is the traslation subgroup of $\Gamma_i$. Hence, $F$ induces an isomorphism $\varphi:H_1\to H_2$. Let $p_i:\Gamma_i\to H_i$ be the natural projections. Let $\gamma_1\in \Gamma_1$ and suppose $p_1(\gamma_1)=h_1\in H_1$. By the definition of $\varphi$, we have $p_2(F(\gamma_1))=\varphi(h_1)$. If $\alpha:=F|_{\Lambda_1}$ then, for all $\gamma\in \Lambda_1$ $\alpha(h\cdot \gamma)=\alpha(\gamma_1 \gamma \gamma_1^{-1})=F(\gamma_1)F(\gamma)F(\gamma_1)^{-1}=\varphi(h_1)\cdot \alpha(\gamma)$.
\end{obs}

	\begin{teo}[Bieberbach's Second Theorem]
		Let $f:\Gamma_1\to \Gamma_2$ be a group isomorphism of two Bieberbach subgroups of $\Iso(\R^n)$. Then there exists $\alpha\in \mathcal{A}_n:=\GL_n(\R)\ltimes\R^n$ such that $f(\beta)=\alpha\beta\alpha^{-1}$ for all $\beta\in\Gamma_1$. 
	\end{teo} 
	
	\begin{teo}[Bieberbach's Third Theorem]  
		For each $n\in\N$, there are only finitely many isomorphism classes of Bieberbach subgroups of $\Iso(\R^n)$.
	\end{teo}
	
	\medskip
	
	We can study compact flat manifolds from studying Bieberbach groups. Indeed, it is well known that a Riemannian manifold $M$ is flat if and only if $M$ is isometric to $\R^n/\Gamma$ where $\Gamma$ is a subgroup of $\Iso(\R^n)$ that acts properly discontinuously on $\R^n$. In conclusion, Bieberbach subgroups of $\Iso(\R^n)$ are just the fundamental groups of compact flat manifolds.
	
Furthermore, the holonomy group of the Riemannian manifold $\R^n/\Gamma$ can be identified with the finite group $\Gamma/\Lambda$ (see for instance \cite[page 50]{Ch}). We will denote $\mathrm{Hol}(\Gamma):=\Gamma/\Lambda$.
	
	As a consequence of Bieberbach's Second and Third we have 
	\begin{teo}
		Let $M$ and $N$ be compact flat manifolds with isomorphic fundamental groups. Then there exists an affine equivalence\footnote{An \textit{affine equivalence} between two Riemannian manifolds $M$ and $N$  is a diffeomorphism $F:M\to N$ such that $f^*\nabla^M =\nabla^N$, where $\nabla^M$ (respectively $\nabla^N$) is the Levi-Civita connection on $M$ (respectively $N$).} between $M$ and $N$. Moreover, for each $n\in\N$ there are only finitely many classes of affine equivalence of compact flat manifolds in dimension $n$.
	\end{teo}
	
	\medskip
	
	We are interested in a special type of compact flat manifolds which arise from solvable Lie groups equipped with a flat left invariant metric. In \cite{Mi}, Milnor gave a nice description of such groups.
	
	\begin{teo}\label{algebraplana}
		A left invariant metric $\la\cdot,\cdot\ra$ on a Lie group $G$ is flat if and only if the associated Lie algebra $\g$ splits as an orthogonal direct sum $\g=\b\oplus\u$, where $\b$ is an abelian subalgebra, $\u$ is an abelian ideal and the linear transformation $\ad_b$ is skew-adjoint for every $b\in\b$.
	\end{teo}
	
	We will call $(G,\la\cdot,\cdot\ra)$ a flat Lie group and  $(\g,\la\cdot,\cdot\ra_e)$ will be called a flat Lie algebra. Using this characterization, Barberis, Dotti and Fino decompose a flat Lie algebra further in the following way \cite[Proposition 2.1]{BDF}. 
	
	\begin{teo}\label{alglieplanas}
		Let $(\g,\la\cdot,\cdot\ra_e)$ be a flat Lie algebra. Then $\g$ splits as an orthogonal direct sum, \[\g=\b\oplus\z(\g)\oplus [\g,\g]\] where $\b$ is an abelian subalgebra, $[\g,\g]$ is abelian and the following conditions are satisfied:
		
		\begin{enumerate}
			\item $\ad:\b\to\mathfrak{so}([\g,\g])$ is injective,
			
			\item $\dim [\g,\g]$ is even, and 
			
			\item $\dim\b\leq\frac{\dim[\g,\g]}{2}$.
			
		\end{enumerate}
	\end{teo}
	
	As a consequence, since $\{\ad_X\mid X\in\b\}$ is an abelian subalgebra of $\so([\g,\g])$, it is contained in a maximal abelian subalgebra. Since these are all conjugate, there exists an orthonormal basis $\B$ of $\z(\g)\oplus[\g,\g]$ and $\lambda_1,\ldots,\lambda_n\in \b^*$ such that for $X\in\b$,  
	\begin{equation}\label{adplana}
		[\ad_X]_\B=\matriz{0_s &&&&&\\&0&-\lambda_1(X)&&&\\&\lambda_1(X)&0&&&\\&&&\ddots&&\\&&&&0&-\lambda_n(X)\\&&&&\lambda_n(X)&0},
	\end{equation}
	where $n=\frac{\dim[\g,\g]}{2}$ and  $s=\dim\z(\g)$.
	
	Note that a flat Lie algebra $(\g,\la\cdot,\cdot\ra_e)$ is $2$-step solvable, since $[\g,\g]$ is abelian, and unimodular\footnote{A Lie algebra $\g$ is said to be \textit{unimodular} if $\tr \ad_X=0$ for all $X\in\g$.}, since $\ad_X$ is skew-adjoint for all $X\in\b$. Also from this it follows that the nilradical of $\g$ is $\z(\g)\oplus [\g,\g]$.
	
	\bigskip
	
	We will look for discrete subgroups of flat Lie groups such that $\Gamma\backslash G$ is compact. This space endowed with the induced flat metric from $G$ is a compact flat manifold. 
	
	In general, if $G$ is a solvable Lie group a discrete and cocompact subgroup $\Gamma$ of $G$ is called a \textit{lattice} and the quotient $\Gamma\backslash G$ is called a \textit{solvmanifold}. With this definition, solvmanifolds are always compact, orientable, and parallelizable.
	
	It is well known that every simply-connected solvable Lie group $G$ is diffeomorphic to $\R^n$ for $n=\dim G$. This implies that the higher homotopy groups of a solvmanifold vanish and $\pi_1(\Gamma\backslash G)\cong\Gamma$.
	
	\begin{obs}
		When $G$ is a flat Lie group and $\Gamma$ is a lattice of $G$ then $\Gamma$ is a Bieberbach group, since $\pi_1(\Gamma\backslash G)=\Gamma$.
	\end{obs}
	
	\medskip
	
	The fundamental group plays an important role in
	the study of solvmanifolds. Indeed, Mostow’s Theorem below shows that solvmanifolds are
	classified, up to homeomorphism, by their fundamental groups. 
	
	\begin{teo}[Mostow] Let $G_1$ and $G_2$ be simply-connected solvable Lie groups with $\Gamma_i$ a lattice in $G_i$ for $i=1, 2$. If $\Phi:\Gamma_1\to \Gamma_2$ is an isomorphism then there exists a diffeomorphism $\tilde{\Phi}: G_1\to G_2$ such that $\tilde{\Phi}|_{\Gamma_1}=\Phi$ and $\tilde{\Phi}(\gamma g)=\Phi(\gamma)\tilde{\Phi}(g)$ for all $\gamma\in\Gamma_1, g\in G_1$. 
	\end{teo} 
	
	\begin{coro}\label{latticedifeo} Two solvmanifolds with isomorphic fundamental groups are diffeomorphic.
	\end{coro}
	
	In particular, two solvmanifolds are diffeomorphic if and only if they are homeomorphic.
	
	\
	
	\section{Splittable flat solvmanifolds}
	
	We are interested in classifying 6-dimensional flat solvmanifolds (up to homeomorphism). However, to determine all the lattices in a given solvable Lie group is a very difficult task and has only been performed for some special cases, for example the $(2n+1)$-Heisenberg group \cite{BG} and the 4-dimensional oscillator group \cite{Fischer}. We will focus on a special type of 6-dimensional flat solvmanifolds, namely the splittable ones. 
	
	Let us begin by recalling some structure theory on solvable Lie groups.
	
	Let $G$ be a simply-connected solvable Lie group, and $N$ the nilradical of $G$ (i.e., the connected closed Lie subgroup of $G$ whose Lie algebra is the nilradical $\n$ of $\g$). Moreover, $[G,G]$ is the connected closed Lie subgroup with Lie algebra $[\g,\g]$. As $G$ is solvable, $[G,G]\subset N$ so $N\backslash G$ is abelian, and from the long exact sequence of homotopy groups associated to the fibration $N\to G\to N\backslash G$ it follows that $N\backslash G$ is simply-connected. Therefore $N\backslash G\cong \R^k$ for some $k\in\N$ and $G$ satisfies the short exact sequence \[1\to N\to G\to \R^k\to 1.\] $G$ is called \textit{splittable} if this sequence splits, that is, there is a right inverse homomorphism of the projection $G\to \R^k$. This condition is equivalent to the existence of a homomorphism $\phi:\R^k\to \Aut(N)$ such that $G$ is isomorphic to the semidirect product $\R^k \ltimes_\phi N$.
	
	Let $\Gamma$ be a lattice in a connected solvable Lie group $G$, and $N$ the nilradical of $G$. Then, the following theorem is well known.
	
	\begin{teo}[Mostow \cite{Mostow1}]
		$\Gamma\cap N$ is a lattice in $N$.
	\end{teo}
	
	\medskip
	
	Following \cite{Y}, a lattice $\Gamma$ of a splittable Lie group $\R^k\ltimes_{\phi} N$ will be called \textit{splittable} if it can be written as $\Gamma=\Gamma_1\ltimes_{\phi} \Gamma_2$ where $\Gamma_1\subset\R^k$ and $\Gamma_2\subset N$ are lattices of $\R^k$ and $N$ respectively. Consequently $\Gamma\backslash G$ will be called a splittable solvmanifold.
	
	\medskip
	
	\begin{obs}
		It is not true that every lattice of a splittable Lie group $\R^k \ltimes_\phi N$ is splittable, as the next example shows. 
	\end{obs} 
	
		\begin{ejem}[A non-splittable lattice in a splittable Lie group]\label{Jonas}
		
		We owe this example to Prof. Jonas Deré. Consider $G=\R^2\ltimes_{\phi} \R^4$. Let $\mathcal{B}=\{e_1,e_2\}$ and\footnote{We will denote $A\oplus B$ the block diagonal matrix $\matriz{A&0\\0&B}$.} \[ [\ad_{e_1}]_{\B'}=\matriz{0&-\pi\\\pi&0}\oplus \matriz{0&0\\0&0},\quad [\ad_{e_2}]_{\B'}=\matriz{0&0\\0&0}\oplus \matriz{0&-\pi\\\pi&0},\] for some basis $\B'$ of $\R^4$. Then \[ \phi(t,s)=\matriz{\cos (\pi t)& -\sin (\pi t) \\ \sin(\pi t)& \cos(\pi t)}\oplus \matriz{\cos (\pi s)& -\sin (\pi s) \\ \sin(\pi s)& \cos(\pi s)}.\] 
		In $G$ consider the subset \[\Gamma=\left\{(z_1, z_2, z_3-\frac{z_2}{2}, z_4, z_5, z_6)\mid z_i\in\Z\right\}.\]
		
		An easy computation shows that for $\gamma_1=(m_1,\ldots,m_6), \gamma_2=(n_1,\ldots,n_6)\in \Gamma$,  
		\begin{align*}
			\gamma_1\gamma_2 &= (m_1+ n_1, m_2+ n_2, m_3\pm n_3 - \frac{m_2\pm n_2}{2}, m_4\pm n_4, m_5 \pm n_5, m_6 \pm n_6)\\
			&= (m_1 + n_1, m_2 + n_2, (m_3\pm n_3 + n_2 \frac{\mp 1+1}{2})-\frac{m_2+n_2}{2}, m_5\pm n_5, m_6\pm n_6)\in \Gamma,\\
			\gamma_1^{-1}&=(-m_1, -m_2, \mp m_3\pm \frac{m_2}{2}, \mp m_4, \mp m_5, \mp m_6)\\
			&=(-m_1, -m_2, \mp m_3 +\frac{\pm m_2-m_2}{2}+ \frac{m_2}{2}, \mp m_4, \mp m_5, \mp m_6)\in \Gamma.
		\end{align*}
		Moreover, it is easily seen that $\Gamma$ is discrete and cocompact, so $\Gamma$ is a lattice of $G$.
		
		If $\Gamma$ were isomorphic to a semidirect product $\Z^2 \ltimes \Z^4$, there would exist elements $\alpha,\beta\in \Gamma$ with $[\alpha,\beta]=e$ such that their projections to $\R^2$ generate $\Z^2$. Since not both $\alpha_1 \beta_2$ and $\alpha_2\beta_1$ can be neither even nor odd (otherwise they would not generate $\Z^2$), we can assume without loss of generality that $\alpha_1\beta_2$ is odd and $\alpha_2\beta_1$ is even. Nevertheless, this contradicts the fact that $[\alpha,\beta]=e_G$.	\end{ejem}
	\medskip

	In \cite{Y}, a useful criterion to determine all the splittable lattices in a splittable Lie group $G=\R^k\ltimes_{\phi} N$ was given. We restate the theorem and for completeness we give a proof in our case of interest, i.e., when $N$ is abelian. 
	
	Note that, fixing a basis $\{X_i\}_{i=1}^k$ of $\R^k$, we have  \begin{equation}\label{expprod}
		\phi( \sum_{i=1}^k t_i X_i)=\prod_{i=1}^k \exp(t_i \ad_{X_i}),\qquad t_i\in\R, \, 1\leq i\leq k.
	\end{equation} Indeed, we have the following commutative diagram:
	\[
	\begin{tikzcd}
		\R^k \arrow{r}{\ad} \arrow[swap]{d}{\mathrm{id}} & \gl_{\dim N}(\R) \arrow{d}{\exp} \\
		\R^k\arrow{r}{\phi} & \GL_{\dim N}(\R)
	\end{tikzcd}
	,\] and $\exp(\ad_{X_i})\exp(\ad_{X_j})=\exp(\ad_{X_j})\exp(\ad_{X_i})$ because $\{\ad_X\mid X\in \R^k\}$ is an abelian subalgebra.
	
	\begin{teo}\label{lattices}
		Let $G=\R^k \ltimes_{\phi} \R^{m}$ be a splittable Lie group, where $\R^m$ is the nilradical of $G$. Then $G$ has a splittable lattice if and only if there exists a basis $\{X_1,\ldots,X_k\}$ of $\R^k$  such that $\exp(\ad_{X_i})$ is similar\footnote{A $n\times n$ matrix $A$ will be said to be \textit{similar} (or \textit{conjugated}) to $B$ if there exists $P\in \GL_n(\R)$ such that $P^{-1}AP=B$ and \textit{integrally similar} if $P\in \GL_n(\Z)$.} to an integer matrix for all $1\leq i\leq k$. In this case, the lattice is $\Gamma=(\bigoplus_{i=1}^k  \Z X_i) \ltimes_{\phi} P \Z^m$ where $P^{-1} \exp(\ad_{X_i}) P$ is an integer matrix.  
	\end{teo}
	\begin{proof}
		$\Leftarrow)$ Let $\Gamma_1$ be the lattice $\bigoplus_{i=1}^k \Z X_i\subset \R^k$ and $\Gamma_2$ the lattice $P\Z^m\subset \R^m$ where $E_i:=P^{-1} \exp(\ad_{X_i})P$ is an integer matrix for all $1\leq i\leq k$. If $\gamma\in \Gamma_1$ then $\gamma=\sum m_i X_i$, $m_i\in \Z$. Therefore, according to \eqref{expprod}, \[\phi(\gamma)\Gamma_2= \left(\prod_{i=1}^k \exp(\ad_{X_i})^{m_i}\right) P \Z^m = P \left(\prod_{i=1}^k E_i^{m_i}\right) \Z^m\subset P \Z^m,\] so $\Gamma_1$ preserves $\Gamma_2$ and therefore the semidirect product $\Gamma=\Gamma_1\ltimes_\phi \Gamma_2$ is well defined. Clearly, $\Gamma$ is a discrete subgroup of $G$. The cocompactness can be seen as in \cite[Theorem 2.4]{Y}. Thus, $\Gamma$ is a lattice of $G$.
		
		$\Rightarrow)$ Let $\Gamma=\Gamma_1\ltimes_{\phi} \Gamma_2$ be a splittable lattice in $G$, where $\Gamma_1\subset \R^k$ and $\Gamma_2\subset \R^m$. Then there exist a basis $\B=\{X_1,\ldots,X_k\}$ and a matrix $P\in \GL_{m}(\R)$ such that $\Gamma_1=\bigoplus_{i=1}^k \Z X_i$ and $\Gamma_2=P\Z^m$. Moreover, as the semidirect product is well defined it follows that for $\gamma\in \Gamma_1$, $\phi(\gamma) P \Z^m\subset P\Z^m$. In particular, choosing $\gamma=X_i$ for $1\leq i\leq k$ we get $\exp(\ad_{X_i}) P \Z^m\subset P \Z^m$. Therefore $P^{-1}\exp(\ad_{X_i})P$ must be an integer matrix for all $1\leq i\leq k$.
	\end{proof}
	
	\smallskip
	
	\begin{obs}\label{isomorfismo}
		A lattice $\Gamma=(\bigoplus_{i=1}^k \Z X_i) \ltimes_{\phi} P \Z^{m}$ as above is isomorphic as a group to $\Sigma_{E_1,\ldots,E_k}:=\Z^k\ltimes_{E_1,\ldots,E_k} \Z^{m}$ where the multiplication is given by \[(r,t)\cdot (r',t')=\left(r+r',t+E_1^{r_1}\cdots E_k^{r_k} t'\right),\; r=(r_1,\ldots,r_k), \;r'\in \Z^k,\, t,t'\in \Z^m.\] Indeed, $f:\Gamma\to \Sigma$ given by $f(\sum r_i X_i, Pt)=(r,t)$ is an isomorphism. The multiplication is well defined because $E_i E_j=E_j E_i$ for all $1\leq i,j\leq k$.
	\end{obs}
	
	\medskip
	
A splittable Lie group $G=\R^k\ltimes_{\phi} \R^m$ with $k=1$ is called an \textit{almost abelian} Lie group in the literature and also can be defined by saying that its Lie algebra $\g$ has a codimension-one abelian ideal. It follows from \cite{Bock} that every lattice in an almost abelian Lie group is splittable. Accordingly, an \textit{almost abelian solvmanifold} is a solvmanifold $\Gamma\backslash G$ such that $G$ is almost abelian.  

\smallskip

\begin{ejem}
In this example we show that there exist solvmanifolds $S=\Gamma'\backslash G'$ with $G'$ non almost abelian such that $S$ is diffeomorphic to an almost abelian solvmanifold, i.e., $\Gamma'$ is isomorphic to a lattice $\Gamma$ in an almost abelian Lie group.

Indeed, let $\g=\R^2\ltimes \R^4$ where \[\ad_{e_1}=\matriz{0&-\pi\\\pi&0}\oplus\matriz{0&-\frac{\pi}{2}\\\frac{\pi}{2}&0},\quad \ad_{e_2}=\matriz{0&-2\pi\\2\pi&0}\oplus\matriz{0&-2\pi\\2\pi&0}.\] There is no codimension-one abelian in $\g$ because $\ad_{e_1}$ and $\ad_{e_2}$ are linearly independent. Therefore $\g$ is not almost abelian. 

Nevertheless, $A=\exp(\ad_{e_1})$ and $B=\exp(\ad_{e_2})$ are integer matrices so $G$ has a splittable lattice $\Gamma$ which is isomorphic to $\Z^2\ltimes_{A,B} \Z^4$. Given that $B=\Id$, the identity map is an isomorphism between $\Z^2\ltimes_{A,B} \Z^4$ and $\Z\ltimes_{(1)\oplus A} \Z^5$, which is isomorphic to a lattice in an almost abelian Lie group. This implies by Corollary \ref{latticedifeo} that $\Gamma\backslash G$ is diffeomorphic to an almost abelian solvmanifold.
\end{ejem}

\medskip	
	
We will observe next that any flat Lie group is a splittable Lie group and our purpose will be to classify splittable lattices in flat Lie groups. 
	
	\smallskip
	
	Let $\g=\b\oplus \z(\g)\oplus [\g,\g]$ be a flat Lie algebra. We may re-write it as $\g=\R^k \ltimes_{\ad} \R^{s+2n}$ where $\b\cong \R^k$ and the nilradical is given by $\z(\g)\oplus [\g,\g]\cong \R^{s+2n}$. Here $s=\dim \z(\g)$ and $2n=\dim [\g,\g]$. Fixing $\{X_1,\ldots,X_k\}$ a basis of $\R^k$, we can write the simply-connected group associated to $\g$ as $G=\R^k\ltimes_\phi \R^{s+2n}$ where $\phi( \sum_{i=1}^k s_i X_i)=\prod_{i=1}^k \exp(s_i \ad_{X_i})$, where $\{\ad_{X_i}\}_{i=1}^k$ are as in \eqref{adplana}. Therefore $\det(\phi(X))=1$ for all $X\in \R^k$.
	
	Note that there may be more than one set $\{E_1,\ldots,E_k\}$ of integer matrices to which we can conjugate $\{\exp(\ad_{X_i})\}_{i=1}^k$. We show next that there is a close relation between splittable lattices of flat Lie groups and finite abelian subgroups of $\GL(s+2n,\Z)$.
	
	\begin{prop}\label{hol}
		Let $G=\R^k\ltimes_\phi\R^{s+2n}$ be a splittable flat Lie group and $\Gamma$ a splittable lattice given by $\Gamma=(\bigoplus_{i=1}^k \Z X_i)\ltimes_{\phi} P \Z^{s+2n}$, where $E_i:= P^{-1}\exp(\ad_{X_i})P$ is integer for $1\leq i\leq k$. Then $\hol(\Gamma\backslash G)\cong \la E_1,\ldots, E_k\ra$.
	\end{prop}
	\begin{proof}
		Recall that $\Gamma\cong \Sigma_{E_1,\ldots,E_k}=\Z^k\ltimes_{E_1,\ldots,E_k}\Z^{s+2n}$. Define $f:\Sigma_{E_1,\ldots,E_k}\to \la E_1,\ldots, E_k\ra$ by $f(r,t)=E_1^{r_1}\cdots E_k^{r_k}$. It is clear that this map is an epimorphism and that \[\Ker f=\{(r_1,\ldots,r_k)\in \Z^k\mid E_1^{r_1}\cdots E_k^{r_k}=\I\}\times \Z^{s+2n}.\] Therefore $\Gamma/\Ker f\cong \la E_1,\ldots, E_k\ra$. 
		To finish, it is enough to prove that $\Ker f$ is the maximal abelian normal subgroup $\Lambda$ of $\Gamma$. It is clear that $\Ker f$ is abelian and normal, so $\Ker f\subset \Lambda$. Conversely, let $(r,t)\in \Lambda$ and $(r',t')\in \Ker f$. As they commute, we have \[(r+r',t+E_1^{r_1}\cdots E_k^{r_k}t')=(r'+r,t'+t).\]
		Varying $t'\in \Z^{s+2n}$ we get $E_1^{r_1} \cdots E_k^{r_k} v=v$ for all $v\in \Z^{s+2n}$, and thus $E_1^{r_1}\cdots E_k^{r_k}=\I$.
	\end{proof}

By Remark \ref{biebiso}, if the Bieberbach groups $\Z^k\ltimes_{E_1,\ldots,E_k} \Z^{s+2n}$ and $\Z^k\ltimes_{F_1,\ldots,F_k} \Z^{s+2n}$ are isomorphic, then their holonomy groups are conjugated in $\GL(k+s+2n,\Z)$. Since they are isomorphic to $\la E_1,\ldots, E_k \ra$ and $\la F_1,\ldots, F_k\ra$ respectively, these subgroups are conjugated in $\GL(s+2n,\Z)$.
The converse statement does not hold, as the next example shows. 

\begin{ejem}
We owe this example to Prof. Derek Holt. Let $A\in \GL(36,\Z)$ be given by

	\[A=\left[\begin{array}{c|c|c}
		v_1 & v_2 &   \begin{array}{c}
			0_{1\times 35} \\
			\hline 
			\I_{35} \\
			\hline 
			w
		\end{array}
	\end{array} \right],\] where \[v_1=-(4,1,4,2,2,4,3,4,1,2,3,3,1,1,4,4,1,2,1,2,2,3,4,4,4,4,2,1,2,4,4,3,2,3,1,-15)^t,\] \[v_2=(149,4,133,64,42,130,76,143,24,53,86,103,35,9,113,144,20,69,22,61,54,82,119,\] \[120,116, 132,68,26,45,118,124,100,47,110,7,120)^t\] and $w=-(1,1,\ldots,1)\in \R^{35}$. This matrix satisfies $A^{37}=\I_{36}$. Using the Magma function \texttt{AreGLConjugate}, it can be seen that $A^2$ is not integrally similar to $A$ or $A^{-1}$. In fact, $A$ is not integrally similar to $A^i$ for all $2\leq i\leq 37$. Moreover, 1 is neither an eigenvalue of $A$ nor of $A^2$. Then, by the following Theorem, $\Sigma_A$ is not isomorphic to $\Sigma_{A^2}$. However $\la A\ra=\la A^2\ra$.
	
	\begin{teo}\cite[Corollary 8.9]{Ventura}
		Let $A,B\in \GL(n,\Z)$ without non-trivial fixed points (i.e. 1 is not an eigenvalue). Then $\Z\ltimes_A \Z^n\cong \Z\ltimes_B \Z^n$ if and only if $B$ is integrally similar to $A$ or $A^{-1}$.
	\end{teo}
	\end{ejem}

Nevertheless, if we forget the matrices and we focus only on the subgroups we have that two conjugate subgroups give rise to isomorphic latticesh, as the following lemma shows.

\begin{lema}
	 Let $E_1,\ldots,E_k,F_1,\ldots,F_k\in \GL(s+2n,\Z)$. If $\la E_1,\ldots, E_k\ra$ is conjugate to $\la F_1,\ldots,F_k\ra$ in $\GL(s+2n,\Z)$ then $\Sigma_{E_1,\ldots,E_k}\cong \Sigma_{F'_1,\ldots,F'_k}$ for some generating set $\{F'_i\}_{i=1}^k$ of $\la F_1,\ldots, F_k\ra$.
\end{lema}
\begin{proof}
	Suppose there exists $Q\in \GL(s+2n,\Z)$ such that $Q^{-1} \la E_1,\ldots, E_k\ra Q= \la F_1,\ldots, F_k\ra$. The matrices $\{Q^{-1} E_i  Q\}_{i=1}^k$ generate $\la F_1,\ldots, F_k\ra$ and $f:\Sigma_{E_1,\ldots,E_k}\to \Sigma_{Q^{-1}E_1 Q,\ldots,Q^{-1}E_k Q}$ given by $f(r,t)=(r,Q^{-1}t)$ is an isomorphism.
\end{proof}

To complete our analysis about the relation between the subgroups $\la E_1,\ldots, E_k\ra$ and the lattices $\Sigma_{E_1,\ldots, E_k}$ we address the case when the cardinality of a minimal generating set of $\la E_1,\ldots, E_k\ra$ is less than $k$. In order to do so we need the following lemma.
	
		\begin{lema}\label{rowreduction}
			Let $E_1,\ldots, E_k\in \GL(m,\Z)$.
		\begin{enumerate}
			\item[\ri] $\Sigma_{E_1,\ldots,E_k}\cong \Sigma_{E_{\sigma(1)},\ldots,E_{\sigma(k)}}$, for all $\sigma\in S_k$. 
			
			\item[\rii] $\Sigma_{E_1,\ldots,E_k}\cong \Sigma_{E_1,\ldots,E_i^{-1},\ldots,E_k}$ for all $1\leq i\leq k$.
			
			\item[\riii] $\Sigma_{E_1,\ldots, E_i,\ldots, E_j,\ldots, E_k}\cong \Sigma_{E_1,\ldots,E_i,\ldots,E_i E_j,\ldots,E_k}$
\end{enumerate} 
	\end{lema}

\begin{proof}
		(i) Given $\sigma\in S_k$, an isomorphism is $f:\Sigma_{E_1,\ldots,E_k}\to\Sigma_{E_{\sigma(1)},\ldots,E_{\sigma(k)}}$ given by $f((r_1,\ldots,r_k),t)=((r_{\sigma(1)},\ldots,r_{\sigma(k)}),t)$.
		
		(ii) An isomorphism $f$ is given by $f((r_1,\ldots,r_i,\ldots,r_k),t)=((r_1,\ldots,-r_i,\ldots,r_k),t)$. 
		
		(iii) $f((r_1,\ldots,r_i,\ldots,r_j,\ldots,r_k),t)=((r_1,\ldots,r_i-r_j,\ldots,r_j,\ldots,r_k),t)$ is an isomorphism between $\Sigma_{E_1,\ldots,E_i,\ldots,E_j,\ldots,E_k}$ and $\Sigma_{E_1,\ldots,E_i,\ldots,E_i E_j,\ldots,E_k}$.
	\end{proof}
Let $A=(a_{ij})_{i,j}\in \GL(k,\Z)$ and define $F_i=E_1^{a_{i1}}\cdots E_k^{a_{ik}}$. 

Since $\Z$ is a Euclidean domain, $A\in \GL(k,\Z)$ if and only if $A$ can be obtained from the identity matrix $\I_k$ by performing a finite sequence of the following elementary row operations:

(i) interchange two rows of $A$;

(ii) multiply a row of $A$ by $\pm 1$;

(iii) for $r\in \Z$ and $i\neq j$, add $r$ times row $j$ to row $i$.

\begin{obs}\label{rowreduction2}
	The previous lemma allows to pass from $\la E_1,\ldots, E_k\ra$ to $\la F_1,\ldots, F_k\ra$ by performing elementary operations on the matrix $A$ and preserving the isomorphism of the corresponding group $\Sigma$. 
\end{obs}

\medskip

Now we are ready to deal with the question related with the minimal generating set. We will use the following theorem.
 
\begin{teo}\cite{Worley}\label{Worley}
	If $(a_1,\ldots,a_m)=1$ then for all integers $n$ the equations \begin{align*}
		&a_1 x_1+\cdots+a_m x_m=n,\\
		&(x_i,x_j)=1,\quad 1\leq i<j\leq m,
	\end{align*} have infinitely many solutions.
\end{teo} 
 
\begin{prop}\label{minimalgeneratingset}
Suppose that the cardinal of a minimal generating set of $\la E_1,\ldots, E_k\ra$ is $\l <k$. Then $\Z^k \ltimes_{E_1,\ldots,E_k} \Z^{s+2n} \cong \Z^\l \ltimes_{H'_1,\ldots,H'_\l} \Z^{s+2n+k-\l}$, where $H'_i=\matriz{\I_{k-\l}&\\& H_i}$ and $\{H_i\}_{i=1}^\l$ is a generating set of $\la E_1,\ldots, E_k\ra$. 	
\end{prop}	
\begin{proof}
	By the structure theorem of finitely generated abelian groups the group $\la E_1,\ldots, E_k\ra$ must be isomorphic to $\Z_{d_1}\times\cdots \times\Z_{d_\l}$. Then there are matrices $F_1,\ldots, F_\l$ of finite order such that $F_1^{m_1}\cdots F_\l^{m_\l}=\I_{s+2n}$ if and only if $(m_1,\ldots,m_\l)\in (\ord F_1)\Z\times\cdots\times (\ord F_\l)\Z$.
	
	
	Let $a_{i1},\ldots, a_{i\l}, b_{j1},\ldots, b_{jk}\in \Z$ such that $E_i=F_1^{a_{i1}}\cdots F_\l^{a_{i\l}}$ and $F_j=E_1^{b_{j1}}\cdots E_k^{b_{jk}}$. Then we have \begin{align*}
		F_1&=(F_1^{a_{11}}\cdots F_\l^{a_{1\l}})^{b_{11}}\cdots (F_1^{a_{k1}}\cdots F_\l^{a_{k\l}})^{b_{1k}}\\
		&= F_1^{\sum_{r=1}^k a_{r1}b_{1r}}\cdots F_j^{\sum_{r=1}^k a_{rj}b_{1r}} \cdots F_k^{\sum_{r=1}^k a_{r\l}b_{1r}}
	\end{align*}  
Therefore there exists $k_1\in \Z$ such that $\sum_{r=1}^k a_{r1}b_{1r}+k_1 \ord F_1=1$. This means that $\operatorname{gcd}(a_{11},\ldots,a_{k1},\ord F_1)=1$. By Theorem \ref{Worley} there are $x_{11},\ldots,x_{k1}, x_1$ pairwise coprime such that $\sum_{r=1}^k a_{r1} x_{r1}+x_1 \ord F_1=1$.

As $\operatorname{gcd}(x_{11},\ldots,x_{k1})=1$, it is known (see for instance \cite[Corollary 3.4.9]{Weintraub}) that there exists a $k\times k$ matrix $A$ in $\SL(k,\Z)$ with first row equal to $(x_{11} \cdots x_{k1})$. Applying Remark \ref{rowreduction2} with $A$  we get $\Sigma_{E_1,\ldots,E_k}\cong \Sigma_{G_1,\ldots,G_k}$ where each $G_i$ is a product of powers of the $F_i's$ but the power of $F_1$ in $G_1$ is equal to 1. Then, by applying Lemma \ref{rowreduction}(iii) we can assume that the power of $F_1$ in all $G_i$ is equal to one. Since the determinant of $A$ is equal to one, $\la G_1,\ldots, G_k\ra=\la E_1,\ldots, E_k\ra=\la F_1,\ldots, F_\l\ra$. We can repeat the process with $F_2$ and obtain $G_1',\ldots,G_k'$ where $G_3',\ldots,G_k'$ do not have powers of $F_1$ and $F_2$. 

Continuing like this up to step $\l$ we will get matrices $G_1^{(\l)}, \ldots, G_{\l-1}^{(\l)}, F_\l, \I_k, \ldots, \I_k$. We rename $H_i:=G_i^{(\l)}, H_\l:=F_\l$ and $H'_i:=\matriz{\I_{k-\l}&\\&H_i}$, then we have $\Sigma_{E_1,\ldots,E_k}\cong \Sigma_{H_1,\ldots,H_\l,\I_k,\ldots,\I_k}$ and it is easy to see that $\Sigma_{H_1,\ldots,H_\l,\I_k,\ldots,\I_k}= \Z^{k-\l}\oplus (\Z^\l\ltimes_{H_1,\ldots,H_\l}\Z^{s+2n})=\Z^\l\ltimes_{H'_1,\ldots,H'_\l} \Z^{s+2n+k-\l}$.
\end{proof}

\begin{coro}
	The holonomy group of a splittable flat solvmanifold $G\backslash \Gamma$ is cyclic if and only if the solvmanifold is diffeomorphic to an almost abelian solvmanifold. 
\end{coro}
	
	
In conclusion, to determine all the diffeomorphism classes of splittable flat solvmanifolds we must determine all the isomorphism classes of splittable lattices. In order to do so, we must look at the finite abelian subgroups of $\GL(n,\Z)$. Two conjugated subgroups give rise to isomorphic lattices, but not all such subgroups are realised as the holonomy group of a flat solvmanifold, as we will see later. However, the problem of classifying the finite abelian subgroups of $\GL(n,\Z)$ up to conjugation for an arbitrary $n$ becomes very difficult. As far as we know, it has only been obtained for $n\leq 6$, as a particular case of the classification of the finite subgroups of $\GL(n,\Z)$  for $n\leq 6$ obtained with the aid of CARAT (see \cite{Pl}). A list of these subgroups can be found in the Internet in https://www.math.kyoto-u.ac.jp/~yamasaki/Algorithm/RatProbAlgTori/crystdat.html.
		
	\medskip

Nevertheless, to classify splittable lattices of flat Lie groups (regardless of the classifications done in low dimensions) one can classify the integral similarity classes of integer matrices, which can be done following the ideas of \cite{Yang}. From this classification the finite abelian subgroups of $\GL(n,\Z)$ can be obtained by doing a careful analysis case by case, although we must be careful because two subgroups may be conjugated without the matrices being conjugated.
	
	\subsection{Almost abelian flat Lie groups}
	In \cite{Tol} we described the structure of an almost abelian flat Lie algebra.
	
	\begin{teo}\label{casiabelianaplana}\cite[Theorem 3.3]{Tol}
		Let $\g=\b\oplus\z(\g)\oplus [\g,\g]$ be a flat Lie algebra. Then $\g$ is almost abelian if and only if $\dim\b=1$. 
	\end{teo}
	
	Then we can write an almost abelian flat Lie algebra $\g$ as $\g=\R x\ltimes_{\ad_x} \R^{s+2n}$ with $s=\dim \z(\g)$, $2n=\dim [\g,\g]$ and in some basis $\B$ of $\R^{s+2n}$ we have \[ [\ad_x]_\B=\matriz{0_s &&&&&\\ & 0&-a_1&&&\\ & a_1&0 &&&\\ &&& \ddots&&\\ &&&& 0&-a_n\\&&&& a_n&0},\quad a_1,\ldots,a_n\in\R\setminus\{0\}.\] The corresponding Lie group $G$ can be written as $G=\R\ltimes_{\phi} \R^{s+2n}$ where \[\phi(t)=\exp(t\ad_x)=\matriz{\I_s&&&\\ & \theta(a_1 t)&&&\\&& \ddots&\\ &&&\theta(a_n t)},\, \text{where}\, \theta(t)=\matriz{\cos t&-\sin t\\ \sin t&\cos t}.\]
	
	We want to determine all the lattices of $G$. According to Theorem \ref{lattices} we have to find $t_0\neq 0$ such that $\phi(t_0)$ is similar to an integer matrix. In view of Proposition \ref{hol}, $\phi(t_0)$ must have finite order.
	
	Computing the integral similarity classes of matrices obtained by conjugating $\phi(t_0)$ (by matrices in $\GL(s+2n,\R)$) actually gives us the integral similarity classes of matrices of finite order in $\SL(s+2n,\Z)$ (recall that $\det \phi(t)=1$ for all $t\in\R$), according to the following theorem.
	
	\begin{teo}\cite{Koo}
		A matrix $A\in\GL(k,\R)$ has finite order if and only if $A$ is similar to \[\I_{k_1}\oplus (-\I_{k_2})\oplus \theta(t_1)^{d_1}\oplus \cdots \oplus \theta(t_r)^{d_r},\] where $k_1,k_2\geq 0$, $r\geq 0$, $d_1,\ldots, d_r\geq 1$, $0<t_1<\cdots<t_r<\pi$, each $t_i$ is a rational multiple of $2\pi$, and $k_1+k_2+2(d_1+\cdots+d_r)=k$.
	\end{teo}
	
	Denote $\tilde{\phi}(t)=\theta(a_1 t)\oplus\cdots\oplus \theta(a_n t)$. We are going to describe a method to find all the possible values of the set $\{a_i t_0\}_{i=1}^n$ such that $\tilde{\phi}(t_0)$ is similar to an integer matrix in terms of the possible integer characteristic polynomials:
	
	\bigskip
	
	\centerline{\texttt{Method to find $\{a_i t_0\}_{i=1}^n$}}
	
	\medskip
	 
	If $n=1$, by looking at the trace of $\phi(t_0)$ it is easy to see that $\phi(t_0)$ is similar to an integer matrix if and only if $a_1 t_0\in\{\frac{\pi}{2}, \frac{3\pi}{2}, \frac{2\pi}{3}, \frac{4\pi}{3}, \frac{\pi}{3}, \frac{5\pi}{3}, \pi, 2\pi\}+2\pi\Z$.
	
	If $n>1$ note that the characteristic polynomial\footnote{Given a matrix $A$, $P_A$ and $M_A$ will denote the characteristic and the minimal polynomials of $A$ respectively.}  $P_{\phi(t_0)}$ is equal to $(x-1)^s P_{\tilde{\phi}(t_0)}$ and thus $P_{\phi(t_0)}\in \Z[x]\iff P_{\tilde{\phi}(t_0)}\in \Z[x]$. Also, if we denote $\bar{\phi}(t_0)=-\I_r\oplus \tilde{\phi}(t_0)$ ($r\in \N$) then we have $P_{\bar{\phi}(t_0)}\in \Z[x]\iff P_{\tilde{\phi}(t_0)}\in \Z[x]$. Therefore we can work with $P_{\tilde{\phi}(t_0)}$ and assume $a_i t_0\notin \{\pi, 2\pi\}+2\pi\Z$ for all $i$. 
	
	Now, as $\tilde{\phi}(t_0)$ has finite order, we have to analyze the integer polynomials $p(x)$ of degree $2n$ such that $p(x)$ divides $x^d-1$ for some $d\in \N$ and $p(x)$ has no real roots. Equivalently, we look for the polynomials of degree $2n$ which can be written as a product of cyclotomic polynomials of degree $\geq 2$. 
	Given that the degree of the $j$-cyclotomic polynomial $\Phi_{j}$ is $\varphi(j)$, where $\varphi$ is the Euler's totient function, we must determine the sets with repetition (also called multisets) $S\subset\{3,4,\ldots\}$ which satisfy \[\sum_{j\in S}\varphi(j)=2n.\] 
	It is known that the totient function has lower bounds such as $\varphi(j)\geq \sqrt{\frac{j}{2}}$ so we can assume that $S\subset \{3,4,\ldots,8n^2\}$, where the number of elements in $S$ is at most $n$. Moreover, the list of integers  $x$ which solve the equation $\varphi(x)=j$ for $1\leq j\leq 1000$ is known\footnote{See http://primefan.tripod.com/TotientAnswers1000.html}.
	
	Assume the matrix $\theta_{i,k}(t_0):=\theta(a_i t_0)\oplus \cdots\oplus \theta(a_k t_0)$ has $P_{\theta_{i,k}(t_0)}=\Phi_j$ for some $j$ and $1\leq i\leq k\leq n$. Then we can deduce the values of $\{a_i t_0,\ldots,a_k t_0\}$. Indeed, looking at the eigenvalues of $\theta_{i,k}(t_0)$ we have \[\left\{\exp\left(\pm \sqrt{-1} a_m t_0\right) \right\}_{m=i}^k=\left\{\exp\left(\frac{2\pi i }{j}\l\right)\mid (\l,j)=1\right\}.\] 
	Then, it is easy to verify that \[\{\exp( \sqrt{-1} a_i t_0)\}_{m=i}^k=\begin{cases} \{\exp(\frac{2\pi i}{j} \l)\mid 1\leq \l \leq \frac{j-1}{2}\},& j\,\text{odd}, \\ \{\exp(\frac{2\pi i}{j} \l)\mid 1\leq \l\leq\frac{j-2}{2}\},& j\,\text{even}\end{cases}.\]
	so $a_i t_0\in \frac{2\pi}{j}\{\l,j-\l\}+2\pi\Z$ for some $\l$ as above.

	\medskip
	
	\section{Splittable flat solvmanifolds of dimension 6}
	
	Our goal is to classify splittable flat solvmanifolds of dimension 6. 
	
	Let $\g=\b\oplus\z(\g)\oplus[\g,\g]$ be a non abelian flat Lie algebra of dimension 6. There are two possibilities for $\dim \b$, namely $\dim \b=1$ or $\dim \b=2$. If $\dim\b=1$, i.e. $\g$ is almost abelian, then $\dim[\g,\g]$ can be 2 or 4 and if $\dim\b=2$ then $\dim[\g,\g]=4$.
	
	\subsection{The almost abelian case $\R\ltimes \R^5$}
	
	A 6-dimensional almost abelian Lie algebra $\g$ can be written as $\g=\R x\ltimes_{\ad_x}\R^5$ where $\ad_x$ can be written in some basis $\mathcal{B}$ of $\z(\g)\oplus [\g,\g]$ as the block matrix \[[\ad_x]=\matriz{0&&&&\\& 0&-a&&\\& a&0&&\\ &&& 0&-b\\&&& b&0},\quad a\in \R,\; b\in \R\setminus\{0\}.\]
	
	The corresponding Lie group is 
	$G=\R\ltimes_{\phi} \R^5$ with $\phi(t)=\I_1\oplus \theta(a t)\oplus \theta(bt)$.
	
	We want to find all the lattices of $G$. In order to do so we have to:
	\begin{enumerate}[(I)]
		\item Find the values $t_0\neq 0$ such that $\phi(t_0)$ is similar to an integer matrix (Theorem \ref{lattices}).
		
		\item For each of the values of $t_0$ found in (1) we have to determine\footnote{We include the computations we made to obtain the classification of the integral similarity classes of integer matrices obtained by conjugating $\phi(t_0)$ because we believe that these calculations are useful for tackling the same problem in higher dimensions where there are no classifications available.} the integral similarity classes of integer matrices obtained by conjugating $\phi(t_0)$. These classes will give us actually the integral similarity classes of matrices with finite order in $\SL(5,\Z)$.
		
		\item For the integral similarity classes we obtained in (II) we have to check if the natural map between the set of these classes and the set of conjugacy classes of finite cyclic subgroups is a bijection, i.e. whether we obtained or not non integrally similar matrices $A,B\in \SL(5,\Z)$ such that $\la A\ra$ and $\la B\ra$ are $\GL(5,\Z)$-conjugate.
	\end{enumerate}
	
	\
	
	(I) \textit{Determining the values of $t_0$}: it can be done as in \cite[Lemma 5.5]{Tol} or applying the method described at the end of the previous section. We obtain the following values.
	
	\begin{prop}\label{lattices5}
		Let $G=\R\ltimes_\phi \R^5$. Then $\phi(t_0)\neq \I_5$ is similar to an integer matrix if and only if one of the following cases occurs:
		
		\begin{center}
			\begin{tabular}{|c|c|c|}
				\hline
				& Values for $at_0$ & Values for $bt_0$\\\hline
				\rule{0pt}{3ex} {Case $(1)$}&$2\pi\Z$& $\{\pi, \frac{2\pi}{3},\frac{4\pi}{3},\frac{\pi}{2},\frac{3\pi}{2},\frac{\pi}{3},\frac{5\pi}{3}\}+2\pi\Z$\\[1ex] \hline
			\end{tabular}
			
			\begin{tabular}{|c|c|c|}\hline
				& Values for $at_0$ & Values for $bt_0$\\\hline
				\rule{0pt}{3ex} {Case $(2)$}& $\{\pi\}+2\pi \Z$&  $\{\pi, \frac{2\pi}{3},\frac{4\pi}{3},\frac{\pi}{2},\frac{3\pi}{2},\frac{\pi}{3},\frac{5\pi}{3}\}+2\pi\Z$\\ [1ex] \hline
			\end{tabular}
			
			\begin{tabular}{|c|c|c|}\hline
				& Values for $at_0$ & Values for $bt_0$\\\hline
				\rule{0pt}{3ex}{Case $(3)$} &$\{\frac{2\pi}{3},\frac{4\pi}{3},\frac{\pi}{2},\frac{3\pi}{2},\frac{\pi}{3},\frac{5\pi}{3}\}+2\pi\Z$ &  $\{\frac{2\pi}{3},\frac{4\pi}{3},\frac{\pi}{2},\frac{3\pi}{2},\frac{\pi}{3},\frac{5\pi}{3}\}+2\pi\Z$  \\[1ex]\hline
			\end{tabular}
			
			\begin{tabular}{|c|c|c|}\hline
				& Values for $at_0$ & Values for $bt_0$\\\hline
				\rule{0pt}{3ex} \multirow{4}{*}{Case $(4)$}
				&$\{\frac{\pi}{4},\frac{7\pi}{4}\}+2\pi\Z$&$\{\frac{3\pi}{4},\frac{5\pi}{4}\}+2\pi\Z$\\\rule{0pt}{3ex}
				&$\{\frac{\pi}{5},\frac{9\pi}{5}\}+2\pi\Z$&$\{\frac{3\pi}{5},\frac{7\pi}{5}\}+2\pi\Z$\\\rule{0pt}{3ex}
				&$\{\frac{2\pi}{5},\frac{8\pi}{5}\}+2\pi\Z$&$\{\frac{4\pi}{5},\frac{6\pi}{5}\}+2\pi\Z$\\\rule{0pt}{3ex}
				&$\{\frac{\pi}{6},\frac{11\pi}{6}\}+2\pi\Z$&$\{\frac{5\pi}{6},\frac{7\pi}{6}\}+2\pi\Z$\\ [1ex]
				\hline
			\end{tabular}	
		\end{center}
	\end{prop}
	\begin{obs}\label{casos} Note that if we change $a t_0$ for $a t_0+2\pi k$ where $k\in \Z$ we do not change the lattice because the matrix $\phi(t_0)$ is the same, and if we change $a t_0$ for $2\pi - at_0$ or the order of the blocks  we get a conjugate matrix to $\phi(t_0)$ (in the first case because $\theta(2\pi -at_0)=\theta(at_0)^{-1})$. As we are looking for \textit{all} the integer matrices for which we can conjugate the matrix $\phi(t_0)$ we can work only with the following matrices $\phi(t_0)$: 
		\begin{flalign*}
			\text{Case (1):} &\quad  \I_3\oplus -\I_2, \quad \I_3\oplus\,\theta\left(\frac{2\pi}{3}\right),\quad \I_3\oplus\,\theta\left(\frac{\pi}{2}\right),\quad \I_3\oplus\,\theta\left(\frac{\pi}{3}\right),\\
			\text{Case (2):} &\quad \I_1\oplus-\I_4,\quad \I_1\oplus-\I_2\oplus\, \theta\left(\frac{2\pi}{3}\right),\quad \I_1\oplus-\I_2\oplus\,\theta\left(\frac{\pi}{2}\right),\quad\I_1\oplus-\I_2\oplus\,\theta\left(\frac{\pi}{3}\right), \\
			\text{Case (3):}&\quad \I_1\oplus\,\theta\left(\frac{2\pi}{3}\right)\oplus\theta\left(\frac{2\pi}{3}\right),\quad \I_1\oplus\,\theta\left(\frac{2\pi}{3}\right)\oplus\theta\left(\frac{\pi}{2}\right),\quad \I_1\oplus\,\theta\left(\frac{2\pi}{3}\right)\oplus\theta\left(\frac{\pi}{3}\right),\\
			&\quad\I_1\oplus\,\theta\left(\frac{\pi}{2}\right)\oplus\,\theta\left(\frac{\pi}{2}\right),\quad 
			\I_1\oplus\,\theta\left(\frac{\pi}{2}\right)\oplus\theta\left(\frac{\pi}{3}\right),\quad
			\I_1\oplus\,\theta\left(\frac{\pi}{3}\right)\oplus\theta\left(\frac{\pi}{3}\right),\\
			\text{Case (4):}&
			\quad \I_1\oplus\, \theta\left(\frac{\pi}{4}\right)\oplus\theta\left(\frac{3\pi}{4}\right),\quad \I_1\oplus\,
			\theta\left(\frac{\pi}{5}\right)\oplus\theta\left(\frac{3\pi}{5}\right),\\
			&\quad\I_1\oplus\,
			\theta\left(\frac{2\pi}{5}\right)\oplus\theta\left(\frac{4\pi}{5}\right),\;\I_1\oplus\,\theta\left(\frac{\pi}{6}\right)\oplus\theta\left(\frac{5\pi}{6}\right).
		\end{flalign*}
	\end{obs}
	
	\
	
	(II)  \textit{Finding the integral similarity classes of integer matrices to which $\phi(t_0)$ is similar}.
	
	\
	
	$\bullet$ \textit{Case} $\ord\phi(t_0)=2$. From the matrices we obtained above, the unique ones with order 2 are $\I_3\oplus -\I_2$ and $\I_1\oplus -\I_4$.
	
	The problem of determining the integral conjugacy classes of integer matrices $A$ such that $A^2=\I$ was solved in complete generality by Hua and Reiner \cite{Hua}.
	
	\begin{lema}
		Every matrix $A\in M_n(\Z)$ such that $A^2=\I_n$ is integrally similar to a matrix of the form \[W(x,y,z)=\underbrace{L\oplus\cdots\oplus L}_{x}\oplus(-\I_y)\oplus\I_z,\] where $2x+y+z=n$ and $L=\matriz{1&0\\1&-1}$.
	\end{lema}
	In our case $n=5$ and since $\det\phi(t_0)=1$ and $\det L=-1$ it follows $x\equiv y\,(2)$. Then the possible values are $(x,y,z)=(0,0,5)$, (0,4,1), (0,2,3), (1,3,0), (1,1,2) or (2,0,1). 
	
	Thus, a set of representatives for integral conjugacy classes of matrices of order 2 in $\SL_5(\Z)$ is \begin{equation}\label{mat5}
		-\I_4\oplus \I_1,\quad -\I_2\oplus \I_3,\quad L\oplus -\I_3,\quad L\oplus -\I_1\oplus \I_2,\quad L\oplus L\oplus \I_1.\end{equation}
	The matrix $\I_3\oplus -\I_2$ is obviously conjugated to $-\I_2\oplus \I_3$ but also to $L\oplus -\I_1\oplus \I_2$ and $L\oplus L\oplus \I_1$ because they have the same Jordan form. The matrix $\I_1\oplus-\I_4$ is conjugated to $-\I_4\oplus\I_1$ and to $L\oplus-\I_3$ because they have the same Jordan form. Therefore with $\phi(t_0)=\I_3\oplus -\I_2$ and $\phi(t_0)=\I_1\oplus-\I_4$ we can conjugate to any of the matrices in \eqref{mat5}.
	
	
	\
	
	$\bullet$ \textit{Case} $\ord\phi(t_0)\neq 2$. We introduce some tools developed by Yang in \cite{Yang}. 
	
	\begin{lema}\label{blocktriangular}\cite[Lemma 1.1]{Yang}
		Every matrix $A\in M_n(\Z)$ is integrally similar to a block triangular matrix \[\matriz{A_{11}&A_{12}&\cdots&A_{1r}\\0&A_{22}&\cdots&A_{2r}\\\vdots&\vdots&\ddots&\vdots\\0&0&\cdots&A_{rr}},\] where the characteristic polynomial of $A_{ii}$ is irreducible for $1\leq i\leq r$ (then $r$ is the number of distinct irreducible factors of $P_A$). The block triangularization can be attained with the diagonal blocks in any prescribed order.
	\end{lema}
	Given a monic polynomial $p(x)\in \Z[x]$ of degree $n$ such that $p(0)=\pm 1$, let $\mathcal{M}_p$ be the set of integral conjugacy classes of matrices $A$ with $P_A=p$, and $\M_A=\{B\in M_n(\Z): P_B(x)=P_A(x),\; M_B(x)=M_A(x)\}$. It is known that $C_p$, the companion matrix of $p(x)$, has characteristic polynomial equal to $p$ so $\M_p\neq\emptyset$. Let $|\M_p|$ be the cardinal of $\M_p$. 
	
	\begin{lema}\label{potciclotomica}\cite[Lemma 1.2]{Yang}
		Let $A\in \GL_n(\Z)$ have irreducible minimal polynomial $f(x)$ with $|\mathcal{M}_f|=1$. Then $A$ is integrally similar to $C_f\oplus \cdots\oplus C_f$. That is $|\mathcal{M}_{f^k}|=1$.
	\end{lema}
	
	A particular case is when $f=\Phi_n$ for certain values of $n$. According to Latimer, MacDuffee and Taussky \cite{LatMac}, if $n<23$ then $|\mathcal{M}_{\Phi_n}|=1$ because $\Q(\xi_n)$, where $\xi_n$ is a primitive $n$-th root of unity, has class number one for $n<23$.
	
	For $A\in \GL_n(\Z)$, we denote $A^{\pm}:=\matriz{A&e_1^\top\\0&\pm 1}$, where $e_1=(1,\ldots,0)$. Let $C_n$ be the companion matrix of $\Phi_n$.
	\begin{teo}\label{caso1}\cite[Theorem 1.4]{Yang}
		Let $n\geq 2$ and $A=\underbrace{C_n\oplus \ldots\oplus C_n}_{s}$. Let $M=\matriz{A&X\\0&\I_m}$, for some integer matrix $X$ of the appropiate size.
		
		\begin{enumerate}
			\item If $n=p^k$, where $p$ is a prime number and $k\geq 1$, then $M$ is integrally similar to  \[\underbrace{C_n^+\oplus\cdots\oplus C_n^+}_{t}\oplus\underbrace{C_n\oplus\cdots\oplus C_n}_{s-t}\oplus \I_{m-t}\] where $t$ satisfies $0\leq t\leq \min(s,m)$ and is uniquely determined by $M$.
			
			\item If $n$ is not a power of a prime, then $M$ is integrally similar to $A\oplus\I_m$.
		\end{enumerate}
	\end{teo}
	
	We analyze the matrices with $\ord\phi(t_0)\neq 2$ obtained in Remark \ref{casos}.
	
	\
	
	\textit{Case (1):} We obtained $\phi(t_0)=\I_3\oplus\, \theta(b t_0)$ where $b t_0\in \{\frac{2\pi}{3}, \frac{\pi}{2}, \frac{\pi}{3}\}$.
	
	We have $P_{\phi(t_0)}(x)=(x-1)^3 \Phi_i(x)$ and $M_{\phi(t_0)}(x)=(x-1) \Phi_i (x)$ for some $i\in \{3,4,6\}$, according to the value of $bt_0$.
	
	By using Lemma \ref{blocktriangular} we obtain that if $B\in \M_{\phi(t_0)}$ then $B$ is integrally similar to a matrix of the form $\matriz{C&X\\0&D}$ where $P_C=\Phi_i$ and $D=\matriz{1&x&y\\0&1&z\\0&0&1}$. Since $M_{\phi(t_0)}(x)=(x-1) \Phi_i(x)$, $D$ must be $\I_3$. Moreover, since $|\M_{\Phi_i}|=1$ for $i\in \{3,4,6\}$ we can assume that $C=C_i$. 
	
	Using Theorem \ref{caso1} (in this case $s=1$ and $m=3$) we get the following table: 
	\begin{center}
		\begin{tabular}{|c|c|}
			\hline
			$i$ & $B$ is integrally similar to\\\hline
			\rule{0pt}{3ex}
			3 & $C_3\oplus \I_3$, $C_3^+ \oplus \I_2$ \\[1ex]\hline
			\rule{0pt}{3ex}
			4 & $C_4\oplus \I_3$, $C_4^+\oplus \I_2$ \\[1ex]\hline
			\rule{0pt}{3ex}
			6 & $C_6\oplus \I_3$\\  \hline
		\end{tabular}
	\end{center}
	
	As the Jordan forms of $A^{\pm}$ and $A\oplus (\pm \I_1)$ are equal it follows that in each case $\phi(t_0)$ is similar to the corresponding matrices in the right column of the table above.
	
	\
	
	\textit{Case (4):} We obtained $\phi(t_0)$ where $(at_0, bt_0)\in \{ (\frac{\pi}{4},\frac{3\pi}{4}), (\frac{\pi}{5},\frac{3\pi}{5}), ( \frac{2\pi}{5},\frac{4\pi}{5}), (\frac{\pi}{6},\frac{5\pi}{6})\}$.
	
	We get $P_{\phi(t_0)}(x)=M_{\phi(t_0)}(x)=(x-1)\Phi_i(x)$ for some $i\in\{5,8,10,12\}$, according to the values of $at_0$ and $bt_0$. By Lemma \ref{blocktriangular}, if $B\in \M_{\phi(t_0)}$ then $B$ is integrally similar to a matrix of the form $\matriz{C&X\\0&\I_1}$. Given that $|\M_{\Phi_i}|=1$ for $i\in\{5,8,10,12\}$ we can assume that $C=C_i$. Using Theorem \ref{caso1}, with $m=s=1$, we get the following table:
	\begin{center}
		\begin{tabular}{|c|c|}
			\hline
			$i$ & $B$ is integrally similar to\\\hline
			\rule{0pt}{3ex}
			5 & $C_5\oplus \I_1 $, $C_5^+$ \\[1ex]\hline
			\rule{0pt}{3ex}
			8 & $C_8\oplus \I_1$, $C_8^+$ \\[1ex]\hline
			\rule{0pt}{3ex}
			10 & $C_{10}\oplus \I_1$ \\[1ex]\hline
			\rule{0pt}{3ex}
			12 & $C_{12}\oplus \I_1$ \\\hline
		\end{tabular}
	\end{center}
	As the Jordan forms of $A^{\pm}$ and $A\oplus (\pm \I_1)$ are equal it follows that in each case $\phi(t_0)$ is similar to the corresponding matrices in the right column of the table above.
	
	\
	
	For cases (2) and (3) we describe the notion of $(A,B)$-equivalence (see \cite[Section 2]{Yang}), which is related to the integral similarity problem for upper block triangular matrices of the form $\matriz{C&X\\0&D}$, where $C,D$ have coprime minimal polynomials. 
	
	Let $A\in\GL_m(\Z)$ and $B\in\GL_n(\Z)$ and suppose that their respective characteristic polynomials $f(x)$ and $g(x)$ are coprime. 
	
	We define a $\Z$-module homomorphism $\psi: M_{m\times n}(\Z)\to M_{m\times n}(\Z)$ given by \[\psi(T)=AT-TB.\] It can be seen that in a suitable basis, the matrix of $\psi$ is $A\otimes \I_n-\I_m\otimes B^{\top}$, where $\otimes$ is the Kronecker product of matrices. Then the determinant of $\psi$ is equal to $R(f,g)$, the \textit{resultant} of $f(x)$ and $g(x)$.
	
	We recall that if $p(x)=x^m+\sum_{i=0}^{m-1} a_i x^i$ and $q(x)=x^n+\sum_{i=0}^{n-1} b_i x^i$ are two monic polynomials in $\Z[x]$ then the resultant $R(p,q)$ is the determinant of the Sylvester matrix which is given by \[
	\mathrm{Syl}(p,q)=\left[\phantom{\begin{matrix}a_0\\ \ddots\\a_0\\b_0\\ \ddots\\b_0 \end{matrix}}
	\right.\hspace{-1.5em}
	\begin{matrix}
		1 & a_{m-1} &\cdots& \cdots & a_0 & \\
		& 	\ddots& \ddots& && \ddots & \\
		&& 1 & a_{m-1} & \cdots & \cdots&a_0 \\
		1 & b_{n-1} & \cdots & b_0 & \\
		&\ddots &\ddots & & \ddots & \\
		&& 1 & b_{n-1}& \cdots & b_0
	\end{matrix}
	\hspace{-1.5em}
	\left.\phantom{\begin{matrix}a_0\\ \ddots\\a_0\\b_0\\ \ddots\\b_0 \end{matrix}}\right]\hspace{-1em}
	\begin{tabular}{l}
		$\left.\lefteqn{\phantom{\begin{matrix} a_0\\ \ddots\\ a_0\ \end{matrix}}}\right\}n$\\
		$\left.\lefteqn{\phantom{\begin{matrix} b_0\\ \ddots\\ b_0\ \end{matrix}}} \right\}m$
	\end{tabular}
	\]
	Since $f(x)$ and $g(x)$ are coprime, then $\det\psi=R(f,g)\neq 0$, so $\psi$ is injective. Denote $\la A,B\ra:=\im \psi$. Then, $\Coker\psi$ has order $r:=|R(f,g)|$ where \[\Coker\psi=M_{m\times n}(\Z)/\la A,B\ra.\]
	
	We define an equivalence relation on $M_{m\times n}(\Z)$. In what follows, given $X\in \GL_k(\Z)$, $\mathsf{C}(X)$ will always denote the centralizer of $X$ in $\GL_k(\Z)$.
	
	\begin{defi}\cite[Definition 2.2]{Yang}
		$X,Y\in M_{m\times n}(\Z)$ are said to be $(A,B)$-\textit{equivalent}, denoted by $X \equiv Y\,(A,B)$ (or $X\equiv Y$), if there exist $P\in \mathsf{C}(A)$ and $Q\in \mathsf{C}(B)$ such that $XQ-PY\in \la A,B\ra$. The set of $(A,B)$-equivalence classes is denoted by $\mathcal{S}(A,B)$.
	\end{defi}
	
	It is clear that if $X-Y\in \la A,B\ra$, then $X\equiv Y\,(A,B)$. The converse is not necessarily true. By defining a group action of $\mathsf{C}(A)\times \mathsf{C}(B)$ on $\Coker \psi$ (see \cite[p. 486]{Yang}) it can be seen the important fact that if $r=1$ then $|\mathcal{S}(A,B)|=1$, and if $r>1$ then $1<|\mathcal{S}(A,B)|\leq r$. 
	
	The importance of the $(A,B)$-equivalence is its connection with the integral similarity problem of upper block triangular matrices.
	
	\begin{lema}\label{relacion}\cite[Lemma 3.1]{Yang}
		Let $A\in M_m(\Z)$ and $B\in M_n(\Z)$ with $(M_A,M_B)=1$. Then, $\matriz{A&X\\0&B}$ and $\matriz{A&Y\\0&B}$ are integrally similar if and only if $X\equiv Y\,(A,B)$.
	\end{lema}	
	
	Using Lemmas \ref{blocktriangular} and \ref{relacion} we will solve cases (2) and (3) by computing the set of $(A,B)$-equivalence classes for certain pairs of matrices. We note that if $P_{\phi(t_0)}$ is a product of more than two irreducible polynomials there is more than one way to write $P_{\phi(t_0)}$ as a product of two coprime polynomials and compute the resultant. By Lemma \ref{blocktriangular} the result does not depend on the way we do it. We do it in such a way that the resultant be as small as possible. 
	
	We will also need in some cases the following analogue result to Theorem \ref{caso1}.
	
	\begin{teo}\label{caso2}\cite[Theorem 1.5]{Yang}
		Let $n\geq 3$ and $A=\underbrace{C_n\oplus\cdots\oplus C_n}_{s}$. Let $M=\matriz{A&X\\0&-\I_m}$.
		\begin{enumerate}
			\item If $n=2p^k$, where $p$ is a prime and $k\geq 1$, then $M$ is integrally similar to \[ \underbrace{C_n^-\oplus \cdots\oplus C_n^-}_{t}\oplus\underbrace{C_n\oplus\cdots\oplus C_n}_{s-t}\oplus (-\I_{m-t}),\] where $t$ satisfies $0\leq t\leq \min(s,m)$ and is uniquely determined by $M$.
			
			\item If $n\neq 2p^k$, then $M$ is integrally similar to $A\oplus (-\I_m)$.
		\end{enumerate}
	\end{teo}
	
	\
	
	We analyze next cases (2) and (3) from Remark \ref{casos}.
	
	\
	
	\textit{Case (2):} We obtained $\phi(t_0)=\I_1\oplus -\I_2\oplus\, \theta(bt_0)$ where $bt_0\in \{\frac{2\pi}{3},\frac{\pi}{2},\frac{\pi}{3}\}$.
	
	\
	
	$\spadesuit$ Let $bt_0=\frac{2\pi}{3}$. Then $P_{\phi(t_0)}(x)=(x-1)(x+1)^2 \Phi_3(x)$ and $M_{\phi(t_0)}(x)=(x^2-1)\Phi_3(x)$.
	
	By Lemma \ref{blocktriangular} a matrix $B\in \M_{\phi(t_0)}$ is integrally similar to a matrix $\matriz{C_3&X\\0&D}$, where $D=\matriz{-\I_2&Y\\0&\I_1}$. We have $r=|R(\Phi_3(x),(x-1)(x+1)^2)|=3$. Noting that $-\I_2=C_2\oplus C_2$ and using Theorem \ref{caso1} we get that $D$ can be integrally similar to $-\I_2\oplus \I_1$ or $C_2^+\oplus -\I_1$.
	
	
	\
	
	For $D=-\I_2\oplus \I_1$, $X=\matriz{x_1&x_2&x_3\\x_4&x_5&x_6} \in \la C_3, D\ra$ if and only if there exists an integer matrix $T=\matriz{t_1&t_2&t_3\\t_4&t_5&t_6}$ such that $C_3 T-TD=X$, and \[X=\matriz{t_1-t_4&t_2-t_5&-t_3-t_6\\ t_1&t_2&t_3-2t_6}\iff x_3+x_6\equiv 0\,(3).\] Then, $\Coker\psi=\{[0],[E_{13}], [2E_{13}]\}$ where $E_{ij}$ is the matrix (in $M_{2\times 3}(\Z)$) with a 1 in the place $(i,j)$ and 0 elsewhere. If we choose $P=-\I_2$ and $Q=D$ then we have $2 E_{13}Q-P E_{13}\in \la C_3, D\ra$ so $2 E_{13}\equiv E_{13}$ and thus $\mathcal{S}(C_3,D)=\{[0], [E_{13}]\}$.
	
	\
	
	For $D=C_2^+\oplus -\I_1$, \[X=C_3 T-TD=\matriz{t_1-t_4&-t_1-t_2-t_5&t_3-t_6\\t_1&t_2-t_4-2t_5&t_3} \iff 2x_1+x_2+2x_4+x_5\equiv 0\,(3).\] Then $\Coker\psi=\{[0], [E_{11}], [2E_{11}]\}$. Since $2E_{11}\I_3+\I_2 E_{11}\in \la C_3, D\ra$, we deduce that $\mathcal{S}(C_3,D)=\{[0], [E_{11}]\}$.
	
	\
	
	In conclusion, $B$ can be integrally similar to \[C_3 \oplus -\I_2\oplus \I_1, \quad \matriz{C_3& E_{13}\\ 0& -\I_2\oplus \I_1},\quad C_3\oplus C_2^+\oplus -\I_1, \quad \matriz{C_3& E_{11}\\ 0& C_2^+\oplus -\I_1}\]
	
	\
	
	$\spadesuit$ Let $bt_0=\frac{\pi}{2}$. Then $P_{\phi(t_0)}(x)=(x-1)(x+1)^2 \Phi_4(x)$ and $M_{\phi(t_0)}(x)=(x^2-1)\Phi_4(x)$. By Lemma \ref{blocktriangular} a matrix $B\in\M_{\phi(t_0)}$ is integrally similar to $\matriz{\I_1&X\\0&D}$ where $D=\matriz{C_4&Y\\ 0&-\I_2}$. We have $r=|R(x-1,\Phi_4(x)(x+1)^2)|=8$. By Theorem \ref{caso2}, $D$ can be integrally similar to  $C_4\oplus -\I_2$ or $C_4^-\oplus -\I_1$.
	
	\
	
	For $D=C_4\oplus -\I_2$, \[\matriz{x_1&x_2&x_3&x_4}=T-TD=\matriz{t_1-t_2&t_1+t_2&2t_3&2t_4} \iff x_1+x_2\equiv x_3\equiv x_4\equiv 0\, (2).\] Then $\Coker\psi=\{[0], [e_1], [e_3], [e_4], [e_1+e_3], [e_1+e_4], [e_3+e_4], [e_1+e_3+e_4]\}$. 
	
	Choosing $P=\I_1$ and $Q=-\I_2\oplus C_4\in \mathsf{C}(C_4\oplus-\I_2)$ we have $e_3 \equiv   e_4$ and $e_1+e_3 \equiv e_1+e_4$. Choosing $P=\I_1$ and $Q=\I_2\oplus\matriz{1&1\\1&0}$ we have $e_3 \equiv e_3+e_4$ and $e_1+e_3 \equiv e_1+e_3+e_4$.
	
	Using that $\mathsf{C}(D)=\left\{\matriz{p&q\\-q&p}\oplus \matriz{p_1&p_2\\p_3&p_4}\right\}\cap \GL_4(\Z)$ it can be seen that there are no more $(\I_1,D)$-equivalences. For example, if $[e_1]=[e_3]$ then there are $Q\in \mathsf{C}(D)$ and $P=\pm \I_1$ such that $e_1 Q-P e_3=\matriz{p&q&\mp 1&0}\in \la \I_1,D\ra$ but $\mp 1\not\equiv 0\,(2)$.  Therefore, $\mathcal{S}(\I_1,D)=\{[0], [e_1], [e_3], [e_1+e_3]\}$.
	
	\
	
	For $D=C_4^-\oplus -\I_1$, \[X=\matriz{t_1-t_2& t_1+t_2& -t_1+2t_3&2t_4 }\iff 
	x_1+x_2+2x_3\equiv 0\, (4),\,  x_4\equiv 0\, (2).\] Then $\Coker\psi=\{[0], [e_1], [e_3], [e_4], [e_1+e_3], [e_1+e_4], [e_3+e_4], [e_1+e_3+e_4]\}$. 
	
	Now, $\mathsf{C}(D)=\left\{\matriz{q_1&-q_2&q_3&q_4\\q_2&q_1&q_2-q_3&-q_4\\0&0&q_1+q_2-2q_3&-2q_4\\0&0&q_5&q_6}\right\}\cap \GL_4(\Z)$. Let $P=\I_1$.
	
	If we choose $q_1=q_3=q_6=1$ and $q_2=q_4=q_5=0$ then we get $e_1\equiv e_1+e_3$ and $e_1+e_4\equiv e_1+e_3+e_4$. Choosing $
	q_1=q_4=q_6=1$ and $q_2=q_3=q_5=0$ we get $e_1\equiv e_1+e_4$. Finally, choosing $
	q_1=q_5=q_6=1$ and $q_2=q_3=q_4=0$ we get $e_4\equiv e_3+e_4$. It can be seen that there are no more $(\I_1,D)$-equivalences so $\mathcal{S}(\I_1,D)=\{[0], [e_1], [e_3], [e_4]\}$.
	
	\
	
	In conclusion $B$ can be integrally similar to 
	\[\I_1\oplus\, C_4\oplus -\I_2,\quad \matriz{1&e_1\\0&C_4\oplus -\I_2},\quad \matriz{1& e_3\\ 0&C_4\oplus -\I_2},\quad \matriz{1& e_1+e_3\\ 0&C_4\oplus -\I_2},\]\[\I_1\oplus\, C_4^-\oplus-\I_1,\quad \matriz{1& e_1\\ 0&C_4^-\oplus-\I_1},\quad \matriz
	{1&e_3\\0&C_4^-\oplus-\I_1},\quad \matriz{1&e_4\\ 0&C_4^-\oplus-\I_1}.\]
	
	\
	
	$\spadesuit$ Let $bt_0=\frac{\pi}{3}$. Then $P_{\phi(t_0)}(x)=(x-1)(x+1)^2 \Phi_6(x)$ and $M_{\phi(t_0)}(x)=(x^2-1)\Phi_6(x)$. By Lemma \ref{blocktriangular} a matrix $B\in \M_{\phi(t_0)}$ is integrally similar to a matrix of the form $\matriz{\I_1&X\\0& D}$, where $D=\matriz{C_6&Y\\0&-\I_2}$. We have $r=|R(x-1,(x+1)^2 \Phi_6(x))|=4$. Moreover, by Theorem \ref{caso2}, $D$ can be integrally similar to $C_6\oplus-\I_2$ or $C_6^-\oplus-\I_1$.
	
	\
	
	For $D=C_6\oplus-\I_2$, \[X=\matriz{t_1-t_2&t_1&2t_3&2t_4}\iff x_3\equiv x_4\equiv 0\,(2).\] 
	Therefore, $\Coker\psi=\{[0],[e_3], [e_4], [e_3+e_4]\}$. With $Q=C_6\oplus C_4$ we get $e_3 Q-e_4\in \la \I_1,D\ra$ and with $Q=C_6\oplus\matriz{1&1\\0&1}$ and $P=\I_1$ we get $e_3 \equiv e_3+e_4$, so $\mathcal{S}(\I_1,D)=\{[0], [e_3]\}$.
	
	\
	
	For $D=C_6^-\oplus-\I_1$, \[X=\matriz{t_1-t_2& t_1& -t_1+2t_3& 2t_4}\iff  x_2+x_3\equiv x_4\equiv 0\,(2).\] We have $\Coker\psi=\{[0],[e_3],[e_4],[e_3+e_4]\}$.
	
	Now, $\mathsf{C}(D)=\left\{\matriz{q_1&-q_2&q_2-2q_3&-2q_4\\q_2&q_1+q_2&q_3&q_4\\0&0&q_1-q_2+3q_3&3q_4\\0&0&q_5&q_6}\right\}\cap \GL_4(\Z)$. Let $P=\I_1$.
	
	Choosing $
	q_1=2, q_2=q_3=0$ and $q_4=q_5=q_6=1
	$ we get $e_3\equiv e_4$ and choosing $q_1=q_4=q_6=1$ and $q_2=q_3=q_5=0$ we get $e_3 \equiv e_3+e_4$, so $\mathcal{S}(\I_1,D)=\{[0], [e_3]\}$.
	
	\
		
	In conclusion, $B$ can be integrally similar to \[\I_1\oplus\, C_6\oplus-\I_2,\quad \matriz{1&e_3\\0&C_6\oplus-\I_2},\quad \I_1\oplus\, C_6^-\oplus-\I_1,\quad \matriz{1&e_3\\0&C_6^-\oplus-\I_1}.\]
	
	\
	
	\textit{Case (3):} We obtained $\phi(t_0)=\I_1\oplus\,\theta(at_0)\oplus\,\theta(bt_0)$ where $at_0,bt_0\in \{\frac{2\pi}{3},\frac{\pi}{2}, \frac{\pi}{3}\}$.
	
	\
	
	$\spadesuit$ Let $at_0=bt_0=\frac{2\pi}{3}$. Then $P_{\phi(t_0)}(x)=(x-1)\Phi_3(x)^2$ and $M_{\phi(t_0)}(x)=(x-1)\Phi_3(x)$. By Lemma \ref{blocktriangular} a matrix $B\in \M_{\phi(t_0)}$ is integrally similar to a matrix of the form $\matriz{D&X\\0&\I_1}$ where $D=\matriz{C_3&Y\\0&C_3}$. By Lemma \ref{potciclotomica}, $D$ is integrally similar to $C_3\oplus C_3$. By Theorem \ref{caso1} with $s=2$ and $m=1$, $B$ can be integrally similar to $C_3\oplus C_3\oplus \I_1$ or $C_3\oplus C_3^+$.
	
	\
	
	$\spadesuit$ Let $at_0=bt_0=\frac{\pi}{2}$. As in the previous $\spadesuit$ we deduce that  $B\in \M_{\phi(t_0)}$ is integrally similar to $C_4\oplus C_4\oplus \I_1$ or $C_4\oplus C_4^+$.
	
	\
	
	$\spadesuit$ Let $at_0=bt_0=\frac{\pi}{3}$. Then  $B\in \M_{\phi(t_0)}$ is integrally similar to $C_6\oplus C_6\oplus \I_1$.
	
	\
	
	$\spadesuit$ Let $at_0=\frac{2\pi}{3}$ and $bt_0=\frac{\pi}{2}$. Then $P_{\phi(t_0)}(x)=M_{\phi(t_0)}(x)=(x-1)\Phi_3(x)\Phi_4(x)$. By Lemma \ref{blocktriangular} a matrix $B\in \M_{\phi(t_0)}$ is integrally similar to a matrix  $\matriz{C_4&X\\ 0&D}$ where $D=\matriz{C_3& Y\\0 &\I_1}$. By Theorem \ref{caso1}, $D$ is integrally similar to $C_3\oplus\I_1$ or $C_3^+$. We have $r=|R(\Phi_4(x),(x-1)\Phi_3(x))|=2$, and since $1<|\mathcal{S}(C_4,D)|\leq r$, we get $|\mathcal{S}(C_4,D)|=2$.
	
	\
	
	For $D=C_3\oplus\I_1$, \[\matriz{x_1&x_2&x_3\\ x_4&x_5&x_6}=\matriz{-t_2-t_4&t_1+t_2-t_5&-t_3-t_6\\t_1-t_5&t_2+t_4+t_5&t_3-t_6}\iff x_3+x_6\equiv 0\, (2).\]
	Therefore $\Coker\psi=\{[0], [E_{13}]\}$ and $\mathcal{S}(C_4,D)=\{[0],[E_{13}]\}$.
	
	\
	
	For $D=C_3^+$, \[\matriz{x_1&x_2&x_3\\x_4&x_5&x_6}=\matriz{-t_2-t_4&t_1+t_2-t_5&-t_1-t_3-t_6\\t_1-t_5&t_2+t_4+t_5&t_3-t_4-t_6}\iff x_2+x_3+x_5+x_6\equiv 0\,(2).\] Then $\Coker\psi=\{[0], [E_{12}]\}$ and $\mathcal{S}(C_4,D)=\{[0], [E_{12}]\}$. Therefore $B$ can be integrally similar to \[ C_4\oplus C_3\oplus \I_1,\quad \matriz{C_4&E_{13}\\0&C_3\oplus \I_1},\quad C_4\oplus C_3^+,\quad \matriz{C_4&E_{12}\\0&C_3^+}.\]
	
	$\spadesuit$ Let $at_0=\frac{2\pi}{3}$ and $bt_0=\frac{\pi}{3}$. Then $P_{\phi(t_0)}(x)=M_{\phi(t_0)}(x)=(x-1)\Phi_3(x)\Phi_6(x)$. By Lemma \ref{blocktriangular} a matrix $B\in \M_{\phi(t_0)}$ is integrally similar to a matrix of the form $\matriz{C_6&X\\0&D}$ where $D=\matriz{C_3&Y\\0&\I_1}$. By Theorem \ref{caso1}, $D$ can be integrally similar to $C_3\oplus \I_1$ or $C_3^+$. We have $r=|R(\Phi_6(x),(x-1)\Phi_3(x))|=4$.
	
	\
	
	For $D=C_3\oplus \I_1$, \[\matriz{x_1&x_2&x_3\\x_4&x_5&x_6}=\matriz{-t_2-t_4&t_1+t_2-t_5&-t_3-t_6\\t_1+t_4-t_5&t_2+t_4+2t_5&t_3}\iff x_1+x_5\equiv x_1+x_2+x_4\equiv 0\,(2).\] Then $\Coker\psi=\{[0], [E_{11}], [E_{12}], [E_{11}+E_{12}]\}$. With $Q=D$ and $P=\I_2$ we have $E_{11}\equiv E_{12}$. Choosing $Q=\I_3$ and $P=C_6$ we have $E_{11}\equiv E_{11}+E_{12}$, so $\mathcal{S}(C_6,D)=\{[0], [E_{11}]\}$.
	
	\
	
	For $D=C_3^+$, \[\matriz{x_1&x_2&x_3\\x_4&x_5&x_6}=\matriz{-t_2-t_4&t_1+t_2-t_5&-t_1-t_3-t_6\\ t_1+t_4-t_5&t_2+t_4+2t_5&t_3-t_4}\iff x_5 \equiv -x_1\equiv x_2+x_4\,(2).\] Therefore $\Coker\psi=\{[0], [E_{11}], [E_{12}], [E_{11}+E_{12}]\}$. If we choose $Q=D$ and $P=\I_2$ we have $E_{11}\equiv E_{12}$ and if we choose $Q=\I_3$ and $P=C_6$ we have $E_{11}\equiv E_{11}+E_{12}$, so $\mathcal{S}(C_6,D)=\{[0], [E_{11}]\}$. 
In conclusion $B$ can be integrally similar to \[C_6\oplus C_3\oplus\I_1,\quad \matriz{C_6&E_{11}\\0&C_3\oplus\I_1},\quad C_6\oplus C_3^+,\quad \matriz{C_6&E_{11}\\0&C_3^+}.\]
	
	$\spadesuit$ Let $at_0=\frac{\pi}{2}$ and $bt_0=\frac{\pi}{3}$. Then $P_{\phi(t_0)}(x)=M_{\phi(t_0)}(x)=(x-1)\Phi_4(x)\Phi_6(x)$. By Lemma \ref{blocktriangular} a matrix $B\in \M_{\phi(t_0)}$ is integrally similar to $\matriz{C_6&X\\0&D}$ where $D=\matriz{C_4&Y\\0&\I_1}$.  By Theorem \ref{caso1}, the matrix $D$ is integrally similar to $C_4\oplus\, \I_1$ or $C_4^+$. Furthermore, $r=|R(\Phi_6(x),\Phi_4(x)(x-1))|=1$ so $B$ can be integrally similar to $C_6\oplus C_4\oplus \I_1$ or $C_6\oplus C_4^+$.
	
	\
	
	(III) For the integral similarity classes of matrices $A\in \SL(5,\Z)$ that we have found, there are no non integrally similar matrices which give rise to conjugated subgroups. Indeed, it can be verified that for the 48 non integrally similar matrices $A$ we have, $A$ is similar to $A^d$ for $(d,\ord A)=1$, so if $\la A\ra$ is conjugated to $\la B\ra$, then $A\sim B^s$ for some $(s,\ord B)=1$ and then $A\sim B$. Therefore there is a bijection between the set of integral similarity classes of integer matrices obtained by conjugating $\phi(t_0)$ and the set of conjugacy classes of finite cyclic subgroups of $\SL(5,\Z)$. These subgroups give rise to 48 6-dimensional almost abelian flat solvmanifolds which are shown in the following table, along with their holonomy group and the abelianization $\Gamma^{ab}:=\Gamma/[\Gamma,\Gamma]$ of the corresponding lattice $\Gamma=t_0 \Z\ltimes P\Z^5$. Recall that we identify $\Gamma$ with $\Sigma_{E}$, where $E=P^{-1}\phi(t_0)P$. The next result shows how to compute $[\Gamma,\Gamma]$ and it is easily obtained. From this, the abelianization can be computed.
	
	\begin{prop}\label{conmutador}
		For $E\in \GL(n,\Z)$, we have $[\Sigma_E,\Sigma_E]=0\Z\oplus \im(\I_n-E)$.
	\end{prop}

Comparing what we obtained with the finite cyclic subgroups of $\SL(5,\Z)$ of the list of all the finite subgroups of $\GL(5,\Z)$ gives us the same non-conjugated subgroups.
	
	
	\begin{table}[H]
		\centering 
		\begin{tabular}{|c|c|c|c|}
			\hline
			$\phi(t_0)$ & Holonomy group & Subgroup & Abelianization of $\Gamma$ \\
			\hline \rule{0pt}{2ex}
			$\I_5$ & $\{e\}$ & $\la \I_5\ra$ & $\Z^6$\\
			\hline \rule{0pt}{2ex}\multirow{3}{*}{
				$\I_3\oplus \theta(\pi)$} &  &  $\la -\I_2\oplus \I_3\ra$ & $\Z^4\oplus \Z_2^2$ \\
			&$\Z_2$ 	& $\la L\oplus-\I_1\oplus \I_2\ra $ & $\Z^4\oplus \Z_2$ \\
			&	& $\la L\oplus L\oplus \I_1\ra$ & $\Z^4$ \\
			\hline \rule{0pt}{2ex}\multirow{2}{*}{$\I_3\oplus \theta(\frac{2\pi}{3})$} & $\Z_3$ & $\la C_3\oplus \I_3\ra$ & $\Z^4\oplus\Z_3$ \\
			& & $\la C_3^+\oplus\I_2\ra$ & $\Z^4$\\
			\hline \rule{0pt}{2ex}\multirow{2}{*}{$\I_3\oplus \theta(\frac{\pi}{2})$} & $\Z_4$ & $\la C_4\oplus\I_3\ra $ & $\Z^4\oplus \Z_2$ \\
			& & $\la C_4^+\oplus\I_2\ra $ & $\Z^4$ \\
			\hline \rule{0pt}{2ex} $\I_3\oplus \theta(\frac{\pi}{3})$ & $\Z_6$ & $\la C_6\oplus\I_3\ra$ & $\Z^4$ \\
			
			\hline \rule{0pt}{3ex}\multirow{2}{*}{$\I_1\oplus \theta(\frac{2\pi}{5})\oplus\theta(\frac{4\pi}{5})$} & $\Z_5$ & $\la C_5\oplus\I_1\ra$ & $\Z^2\oplus\Z_5$ \\
		&& $\la C_5^+\ra$ & $\Z^2$ \\
		\hline \rule{0pt}{3ex}\multirow{2}{*}{$\I_1\oplus \theta(\frac{\pi}{4})\oplus\theta(\frac{3\pi}{4})$} & $\Z_8$ & $\la C_8\oplus\I_1\ra $ & $\Z^2\oplus\Z_2$\\
		&& $\la C_8^+\ra $ & $\Z^2$ \\
		\hline \rule{0pt}{3ex} $\I_1\oplus \theta(\frac{\pi}{5})\oplus\theta(\frac{3\pi}{5})$ & $\Z_{10}$ & $\la C_{10}\oplus\I_1\ra $ & $\Z^2$ \\
		\hline \rule{0pt}{3ex} $\I_1\oplus \theta(\frac{\pi}{6})\oplus\theta(\frac{5\pi}{6})$ & $\Z_{12}$ & $\la C_{12}\oplus\I_1\ra$ & $\Z^2$ \\
		\hline			
	\end{tabular}
	\caption{6-dim. almost abelian flat solvmanifolds (Cases (1) and (4))}
	\label{Case 4}
\end{table}
\begin{table}[H]
\centering
\begin{tabular}{|c|c|c|c|}
			\hline \rule{0pt}{2ex}\multirow{2}{*}{$\I_1\oplus \theta(\pi)\oplus\theta(\pi)$} & $\Z_2$ & $\la\I_1\oplus -\I_4\ra$ & $\Z^2\oplus \Z_2^4$ \\
			&& $\la L\oplus-\I_3\ra$ & $\Z^2\oplus \Z_2^3$ \\
			\hline \rule{0pt}{2ex}\multirow{4}{*}{$\I_1\oplus\theta(\pi)\oplus\theta(\frac{2\pi}{3})$} &  & $\la C_3\oplus -\I_2\oplus \I_1\ra$ & $\Z^2\oplus \Z_2\oplus \Z_6$ \\
			& $\Z_6$& $\la\matriz{C_3&E_{13}\\ 0&-\I_2\oplus\I_1}\ra$ & $\Z^2\oplus \Z_2^2$ \\
			&& $\la C_3\oplus C_2^+\oplus -\I_1\ra$ & $\Z^2\oplus\Z_6$ \\
			&& $\la\matriz{C_3&E_{11}\\ 0&C_2^+\oplus -\I_1}\ra$  & $\Z^2\oplus \Z_2$ \\
			\hline \rule{0pt}{2ex}\multirow{8}{*}{$\I_1\oplus \theta(\pi)\oplus\theta(\frac{\pi}{2})$} &  & $\la\I_1\oplus C_4\oplus-\I_2\ra$ & $\Z^2\oplus\Z_2^3$\\
			&& $\la\matriz{1&e_1\\0&C_4\oplus-\I_2}\ra$  & $\Z^2\oplus\Z_2^2$\\
			&& $\la\matriz{1&e_3\\0&C_4\oplus-\I_2}\ra$ & $\Z^2\oplus\Z_2^2$\\
			&$\Z_4$& $\la\matriz{1&e_1+e_3\\0&C_4\oplus-\I_2}\ra$ & $\Z^2\oplus\Z_2^2$\\
			&& $\la 1\oplus C_4^-\oplus-\I_1\ra$ & $\Z^2\oplus\Z_2\oplus\Z_4$ \\
			&& $\la\matriz{1&e_1\\0&C_4^-\oplus-\I_1}\ra$  & $\Z^2\oplus \Z_2$\\
			&& $\la\matriz{1&e_3\\0&C_4^-\oplus-\I_1}\ra$  & $\Z^2\oplus\Z_2^2$\\
			&& $\la\matriz{1&e_4\\0&C_4^-\oplus-\I_1}\ra$ & $\Z^2\oplus \Z_4$\\
			\hline \rule{0pt}{2ex}\multirow{4}{*}{$\I_1\oplus\theta(\pi)\oplus \theta(\frac{\pi}{3})$} &  & $\la\I_1\oplus C_6\oplus -\I_2\ra$  & $\Z^2\oplus\Z_2^2$ \\
			&$\Z_6$& $\la\matriz{1&e_3\\0&C_6\oplus -\I_2}\ra$  & $\Z^2\oplus\Z_2$\\
			&& $\la\I_1\oplus C_6^-\oplus -\I_1\ra$  & $\Z^2\oplus \Z_2^2$\\
			&& $\la\matriz{1&e_3\\0&C_6^-\oplus-\I_1}\ra$  & $\Z^2\oplus\Z_2$\\
			\hline
			 \rule{0pt}{3ex}\multirow{2}{*}{$\I_1\oplus \theta(\frac{2\pi}{3})\oplus \theta(\frac{2\pi}{3})$} & $\Z_3$ & $\la C_3\oplus C_3\oplus\I_1\ra$ &$\Z^2\oplus\Z_3^2$ \\
			& & $\la C_3\oplus C_3^+\ra$ & $\Z^2\oplus\Z_3$\\
			\hline
			\rule{0pt}{3ex}\multirow{2}{*}{$\I_1\oplus\theta(\frac{\pi}{2})\oplus\theta(\frac{\pi}{2})$} & $\Z_4$ & $\la C_4\oplus C_4\oplus\I_1\ra $ & $\Z^2\oplus\Z_2^2$ \\
			&& $\la C_4\oplus C_4^+\ra$ & $\Z^2\oplus\Z_2$\\
			\hline 
			\rule{0pt}{3ex} $\I_1\oplus\theta(\frac{\pi}{3})\oplus\theta(\frac{\pi}{3})$ & $\Z_6$ & $\la C_6\oplus C_6\oplus\I_1\ra$ & $\Z^2$ \\
			\hline\rule{0pt}{3ex}\multirow{4}{*}{$\I_1\oplus\theta(\frac{2\pi}{3})\oplus\theta(\frac{\pi}{2})$} &  & $\la C_4\oplus C_3\oplus \I_1\ra$ & $\Z^2\oplus\Z_6$ \\
			&$\Z_{12}$& $\la\matriz{C_4&E_{13}\\ 0&C_3\oplus\I_1}\ra$ & $\Z^2\oplus\Z_3$ \\
			&& $\la C_4\oplus C_3^+\ra$  & $\Z^2\oplus\Z_2$ \\
			&& $\la\matriz{C_4&E_{12}\\0&C_3^+}\ra$ & $\Z^2$ \\
			\hline\rule{0pt}{3ex}\multirow{4}{*}{$\I_1\oplus\theta(\frac{2\pi}{3})\oplus\theta(\frac{\pi}{3})$} &  & $\la C_6\oplus C_3\oplus \I_1\ra$  & $\Z^2\oplus\Z_3$ \\
			&$\Z_6 $& $\la\matriz{C_6&E_{11}\\0&C_3\oplus \I_1}\ra$  & $\Z^2\oplus\Z_3$ \\
			&& $\la C_6\oplus C_3^+\ra$ & $\Z^2$\\
			&& $\la\matriz{C_6&E_{11}\\ 0&C_3^+}\ra$  & $\Z^2$\\
			\hline\rule{0pt}{3ex}\multirow{2}{*}{$\I_1\oplus\theta(\frac{\pi}{2})\oplus\theta(\frac{\pi}{3})$} & $\Z_{12}$ & $\la C_6\oplus C_4\oplus \I_1\ra$  & $\Z^2\oplus\Z_2$ \\
			&& $\la C_6\oplus C_4^+\ra$ & $\Z^2$\\
			\hline
			\end{tabular}
			\caption{6-dim. almost abelian flat solvmanifolds (Cases (2) and (3))}
			\label{Case 3}
			
		\end{table}
	
\medskip
	
	\subsection{The non almost abelian case $\R^2\ltimes \R^4$}
	
	Let $\g=\b\oplus\z(\g)\oplus [\g,\g]$ be a non almost abelian flat Lie algebra of dimension 6, i.e., $\z(\g)=0$, $\dim\b=2$ and $\dim [\g,\g]=4$, according to Theorem \ref{alglieplanas} and Theorem \ref{casiabelianaplana}. Then $\g$ can be written as $\g=\R^2\ltimes_{\ad} \R^4$ where $\R^2=\text{span}\{x,y\}$, and in some basis $\B$ of $[\g,\g]$ we have $[\ad_x]=\matriz{0&-a\\a&0}\oplus\matriz{0&-b\\b&0}$ and $[\ad_y]=\matriz{0&-c\\c&0}\oplus\matriz{0&-d\\d&0},$ where $a^2+c^2\neq 0$, $b^2+d^2\neq 0$ and $ad-bc\neq 0$.
	
	 The simply-connected Lie group is $G=\R^2\ltimes_{\phi}\R^4$ where $\phi(tx+sy)=\exp(t\ad_x)\exp(s\ad_y)$.
	
	As stated in Theorem \ref{lattices}, to determine all the splittable lattices in $G$ we have to look for pairs $\{x,y\}$ such that $P^{-1}\exp(\ad_x)P=A$ and $P^{-1}\exp(\ad_y)P=B$ with $A,B\in \GL(4,\Z)$, for some $P\in \GL(4,\R)$. To find such pairs $\{x,y\}$ note that $\exp(\ad_x)=\theta(a)\oplus\theta(b)$ and $\exp(\ad_y)=\theta(c)\oplus\theta(d)$. The values of $a$ and $b$ $($resp. $c$ and $d)$ such that $\theta(a)\oplus \theta(b)$ $($resp. $\theta(c)\oplus\theta(d))$ is conjugate to an integer matrix are as in Proposition \ref{lattices5} and they can be chosen so as to satisfy the conditions above.
	
Keeping in mind that the group $\Sigma_{A,B}=\Z^2\ltimes_{A,B} \Z^4$ is determined by the conjugacy class of the finite abelian subgroup $\la A,B\ra\subset \SL(4,\Z)$, the strategy will be to extract from the list of all the finite subgroups of $\GL(4,\Z)$ the ones which are finite, abelian, 2-generated and contained in $\SL(4,\Z)$. We do not take into account the cyclic ones because they correspond to almost abelian flat solvmanifolds (Proposition \ref{minimalgeneratingset}). Then, we will see which subgroups can be realised as the holonomy group of a splittable flat solvmanifold.
	
	There are 710 finite subgroups of $\GL(4,\Z)$. Only 33 are 2-generated finite abelian subgroups of $\SL(4,\Z)$. 13 of them give rise to groups $\Sigma_{A,B}$ with abelianization of rank 3, which are the following:
	
	$\la \operatorname{diag}(-1,1,-1,1),\operatorname{diag}(1,-1,-1,1)\ra$, $\la \matriz{-1&0\\0&1}\oplus \matriz{0&-1\\-1&0},\matriz{-1&0\\0&1}\oplus\matriz{0&1\\1&0}\ra$, 
	
	$\la\matriz{1&0&0&0\\0&0&0&1\\0&0&-1&0\\0&1&0&0},\matriz{1&0&0&0\\0&-1&-1&0\\0&0&1&0\\0&0&1&-1}\ra$,
	$\la\matriz{1&0\\0&-1}\oplus\matriz{0&-1\\-1&0},(-\I_2)\oplus \I_2\ra$,
	
	$\la\matriz{-1&0&0&0\\0&1&1&0\\0&0&-1&0\\0&0&-1&1}, \matriz{-1&0&0&0\\0&0&0&-1\\0&0&1&0\\0&-1&0&0}\ra$, $\la\matriz{0&0&1&-1\\0&1&0&0\\0&0&-1&0\\-1&0&-1&0},\matriz{1&1&1&0\\0&-1&0&0\\0&0&-1&0\\0&-1&-1&1}\ra$,	
	
	$\la\matriz{1&0&0&0\\0&0&1&-1\\0&0&-1&0\\0&-1&-1&0},\matriz{1&0&0&0\\0&-1&0&0\\0&1&0&1\\0&1&1&0}\ra$, $\la\matriz{-1&0&0&0\\0&1&0&0\\0&-1&0&-1\\0&-1&-1&0},\matriz{-1&0&0&0\\0&0&1&-1\\0&1&0&1\\0&0&0&1}\ra$,
	
	$\la\matriz{0&0&-1&0\\0&0&0&-1\\-1&0&0&0\\0&-1&0&0},\matriz{0&0&1&0\\0&0&0&-1\\1&0&0&0\\0&-1&0&0}\ra$,
	$\la\matriz{0&-1&1&1\\-1&0&-1&-1\\0&0&0&-1\\0&0&-1&0},\matriz{0&-1\\-1&0}\oplus\matriz{0&1\\1&0}\ra$
	
	$\la\matriz{0&1\\1&0}\oplus\matriz{0&-1\\-1&0},\matriz{-1&0&-1&-1\\0&-1&1&1\\0&0&1&0\\0&0&0&1}\ra$, 
$\la\matriz{1&0&2&0\\-1&0&-1&-1\\0&0&-1&0\\-1&-1&-1&0},\matriz{-1&0&0&-2\\1&0&1&1\\1&1&0&1\\0&0&0&1}\ra$,

	$\la\matriz{0&-1\\-1&0}\oplus \matriz{0&-1\\-1&0},\matriz{0&0&-1&0\\0&0&0&-1\\-1&0&0&0\\0&-1&0&0}\ra$.
	
	These groups cannot be realised as the holonomy group of a 6-dimensional flat solvmanifold because Barberis, Dotti and Fino proved in \cite{BDF} that any flat solvmanifold of even dimension is Kähler, so the first Betti number $b_1$ must be even. We obtained the remaining 20 subgroups, which are shown in the following table.
		\begin{table}[H]
		\centering
		\begin{tabular}{|c|c|c|c|}
			\hline \rule{0pt}{2.5ex}
			$\tilde{A},\, \tilde{B}$ & Holonomy group & $\la A,B\ra$ & $\Gamma^{ab}$ \\
			\hline \rule{0pt}{1ex}\multirow{5}{*}{
				$(2\pi, \pi), (\pi, 2\pi)$} &                    & $\la\I_2\oplus-\I_2,\, -\I_4\ra$               & $\Z^2\oplus \Z_2^4$\\
			&  $\Z_2\oplus \Z_2$  & $\la\matriz{\I_2&E_{22}\\0&-\I_2},\, -\I_4\ra$ & $\Z^2\oplus \Z_2^3$\\
			&                     & $\la\matriz{\I_2&\I_2\\0&-\I_2},\, -\I_4\ra$ & $\Z^2\oplus \Z_2^2 $\\
			\hline \rule{0pt}{1.5ex}
			$(2\pi, \pi), (\pi, \frac{2\pi}{3})$ & $\Z_2\oplus \Z_6$ & $\la-\I_2\oplus C_3,\, -\I_4\ra$ & $\Z^2\oplus\Z_2^2 $ \\
			\hline \rule{0pt}{1.5ex}
			$(2\pi,\pi), (\frac{\pi}{2},2\pi)$ & $\Z_2\oplus \Z_4$ & $\la\I_2\oplus C_4,\, -\I_4\ra$ & $\Z^2\oplus\Z_2^3$ \\
			&					& $\la\matriz{\I_2&E_{11}\\0&C_4},\, -\I_4\ra$ & $\Z^2\oplus \Z_2^2$ \\
			\hline \rule{0pt}{1.5ex}
			$(2\pi,\pi), (\frac{\pi}{2},\frac{2\pi}{3})$ & $\Z_4\oplus \Z_6$ & $\la C_3\oplus C_4,\, -\I_4\ra$ & $\Z^2\oplus\Z_2 $ \\
			\hline \rule{0pt}{1.5ex}
			$(2\pi,\pi), (\frac{\pi}{3},\frac{2\pi}{3})$ & $\Z_2\oplus \Z_6$ & $\la C_3\oplus -C_3,\, -\I_4\ra $ & $\Z^2$  \\
			&                   & $\la\matriz{C_3&E_{11}\\0&-C_3},\, -\I_4\ra$ & $\Z^2$ \\
			\hline \rule{0pt}{1.5ex} 
			$(2\pi, \frac{2\pi}{3}), (\frac{2\pi}{3},\pi)$ & $\Z_3\oplus \Z_6$ & $\la\I_2\oplus C_3,\, C_3\oplus -\I_2\ra$ & $\Z^2\oplus \Z_3 $  \\    
			\hline \rule{0pt}{1.5ex}
			$(2\pi,\frac{\pi}{3}), (\frac{\pi}{3},2\pi)$ & $\Z_6\oplus \Z_6$ & $\la\I_2\oplus -C_3, \, C_6\oplus \I_2\ra$ & $\Z^2$ \\
			\hline \rule{0pt}{1.5ex} 
			$(2\pi, \frac{2\pi}{3}), (\frac{2\pi}{3},2\pi)$ & $\Z_3\oplus \Z_3$ & $\la\I_2\oplus C_3,\, C_3\oplus \I_2\ra$ & $\Z^2\oplus \Z_3^2$ \\
			& &  $\la\matriz{\I_2&E_{11}\\0&C_3},\, B_1\ra$ & $\Z^2\oplus \Z_3$\\
			\hline \rule{0pt}{1.5ex}
			$(2\pi, \frac{2\pi}{3}), (\frac{\pi}{3},\pi)$ & $\Z_3\oplus \Z_6$ & $\la\I_2\oplus C_3,\, C_6\oplus -\I_2\ra$ & $\Z^2$  \\
			& & $\la\matriz{\I_2&E_{11}\\0&C_3},\, B_2\ra$ & $\Z^2$\\
			\hline \rule{0pt}{1.5ex}
			$(2\pi, \frac{\pi}{2}), (\frac{\pi}{2},2\pi)$ & $\Z_4\oplus \Z_4$ & $\la\I_2\oplus C_4,\, C_4\oplus \I_2\ra$ & $\Z^2\oplus\Z_2^2$ \\
			& & $\la\matriz{\I_2&E_{11}\\0 &C_4},\, B_3\ra$ & $\Z^2\oplus\Z_2$  \\
			\hline \rule{0pt}{1.5ex}\multirow{5}{*}{
				$(2\pi, \pi), (\frac{\pi}{2},\frac{\pi}{2})$} & & $\la\I_2\oplus -\I_2,\, C_4\oplus C_4\ra$ & $\Z^2\oplus\Z_2^2$ \\
			& $\Z_2\oplus\Z_4$ & $\la\matriz{\I_2&E_{22}\\0&-\I_2},\, B_4\ra$ & $\Z^2\oplus\Z_2^2$ \\
			& & $\la\matriz{\I_2&\I_2\\0&-\I_2},\, C_4\oplus C_4\ra$ & $\Z^2\oplus\Z_2$ \\
			\hline 			
		\end{tabular}
		\caption{6-dimensional splittable non almost abelian flat solvmanifolds}
		\label{Table 3}
	\end{table}
where \[B_1=\matriz{1&1&0&0\\-3&-2&-2&1\\0&0&1&0\\0&0&0&1},\quad B_2=\matriz{-1&1&0&0\\-3&2&-2&1\\0&0&-1&0\\0&0&0&-1},\]\[	B_3=\matriz{1&1&0&0\\-2&-1&-1&1\\0&0&1&0\\0&0&0&1},\quad
B_4=\matriz{1&-2&0&-1\\1&-1&-1&0\\0&0&1&-1\\0&0&2&-1}.\]
	\section{Conclusions}
	We obtained 68 6-dimensional splittable flat solvmanifolds, of which 48 are almost abelian solvmanifolds and 20 are not.
	
	We will discuss the first Betti number $b_1$ of the solvmanifolds we have obtained. Recall that $\Gamma$ is isomorphic to the fundamental group of the solvmanifold $\Gamma\backslash G$. Therefore, by Hurewicz's theorem, $\Gamma^{ab}\cong H_1(\Gamma\backslash G, \Z)$ and thus $b_1(\Gamma\backslash G)$ is equal to the rank of $\Gamma^{ab}$.
	
	By inspection of Table 3, we observe that all 6-dimensional splittable non almost abelian flat solvmanifolds have $b_1$ equal to 2.
	
	On the other hand, the first Betti number of 6-dimensional almost abelian flat solvmanifolds can be computed using \cite[Proposition 4.7]{Bock}: $b_1(\Gamma\backslash G)=6-\mathrm{rank}(\phi(t_0)-\I_5)$. Note that this coincides with the values in Tables 1 and 2.
	
	A long standing problem concerning solvmanifolds is to determine whether the cohomology of a solvmanifold can be computed exclusively using left invariant forms. That is, when there is an isomorphism $H^*(\Gamma\backslash G,\R)\cong H^*(\g)$ where $H^*(\g)$ denotes the Chevalley-Eilenberg cohomology of the Lie algebra.
	
	Regarding this problem, we can obtain from Table 1 examples of solvmanifolds where this isomorphism does not hold. For instance, the Lie group $G=\R \ltimes_{\phi(t_0)} \R^5$ where $at_0=2\pi$ and $bt_0=\frac{2\pi}{3}$ gives the solvmanifold in Table 1 corresponding to $\phi(t_0)=\I_3\oplus \theta(\frac{\pi}{3})$ and has $b_1$ equal to 4 but $\dim H^1(\g)=\operatorname{codim} [\g,\g]=2$. 
	
	\medskip
	
	As a final remark, we would like to point out that the classification of 6-dimensional splittable flat solvmanifolds can be easily generalized to dimension 7, by using the classification of finite subgroups of $\GL(6,\Z)$.

	\
	
	\textbf{Conflict of interest statement}
	
	On behalf of all authors, the corresponding author states that there is no conflict of interest.

	\
	
	\textbf{Data availability}
	Data sharing not applicable to this article as no datasets were generated or analysed during the current study.
	
	\
\end{document}